\documentclass{amsart}

\usepackage{graphicx} 
\usepackage{overpic}

\usepackage{tikz}
\usetikzlibrary{knots}

\theoremstyle{plain}
\newtheorem{thm}{Theorem}[section]
\newtheorem{cor}[thm]{Corollary}
\newtheorem{lemma}[thm]{Lemma}
\newtheorem{prop}[thm]{Proposition}
\newtheorem{claim}{Claim}

\theoremstyle{remark}
\newtheorem{remark}{Remark}

\begin{document}

\title{A presentation of the pure cactus group of degree four}

\author{Takatoshi Hama}
\address{Graduate School of Integrated Basic Sciences, Nihon University, 3-25-40 Sakurajosui, Setagaya-ku, Tokyo 156-8550, Japan}
\email{chta23018@g.nihon-u.ac.jp}

\author{Kazuhiro Ichihara}
\address{Department of Mathematics, 
College of Humanities and Sciences, Nihon University, 3-25-40 Sakurajosui, Setagaya-ku, Tokyo 156-8550, Japan}
\email{ichihara.kazuhiro@nihon-u.ac.jp}

\date{\today}



\begin{abstract}
We give a simple presentation of the pure cactus group $PJ_4$ of degree four.
This presentation is obtained by considering an action of $PJ_4$ on the hyperbolic plane and constructing a Dirichlet polygon for the action.
As a corollary, we provide a direct alternative proof that $PJ_4$ is isomorphic to the fundamental group of the connected sum of five real projective planes.
\end{abstract}

\maketitle

\section{Introduction}

As an analogue of the braid group, the cactus group $J_n$ was introduced in \cite{HENRIQUES-KAMNITZER}, motivated by the study of quantum groups.
We will provide the purely algebraic definition used in this paper in the next section.
It is worth noting that the group itself had been studied earlier in different contexts.
See \cite{DJS03,devadoss2000tessellationsmodulispacesmosaic} for examples.

The cactus group $J_n$ admits a natural surjection onto the symmetric group $S_n$ of degree $n$.
The kernel of this surjection is called the \textit{pure cactus group} $PJ_n$ of degree $n$.
In \cite[Theorem 9]{HENRIQUES-KAMNITZER},  Henriques and Kamnitzer showed that $PJ_n$ is isomorphic to the fundamental group of the Deligne–Mumford compactification $\overline{M_{0,n+1}}(\mathbb{R})$ of the moduli space of real genus-zero curves with $n + 1$ marked points.

In the degree four case, it is known that $\overline{M_{0,5}}(\mathbb{R})$ is homeomorphic to the connected sum of five real projective planes (see \cite[Example 2.5]{EHKR}).
It follows that $PJ_4$ is isomorphic to the fundamental group of this surface, that is,
\[
\left\langle 
\alpha_1, \alpha_2, \alpha_3, \alpha_4, \alpha_5 \mid 
\alpha_1^2 \alpha_2^2 \alpha_3^2 \alpha_4^2 \alpha_5^2 =e
\right\rangle.
\]

On the other hand, another explicit presentation of $PJ_4$ was obtained purely algebraically in \cite[Appendix A]{BCL24}, using the Reidemeister–Schreier method.
This presentation is a simple one-relator presentation with five generators, but it is not related to the presentation described above.
In fact, \cite[Remark 5.6]{BCL24} states that it seems difficult to express the generators $\alpha_k$ explicitly in terms of the (standard) generators $s_{i,j}$ of the full cactus group $J_4$.

In this paper, we provide another presentation of the pure cactus group $PJ_4$ of degree four.
This presentation is obtained geometrically by considering an action of $PJ_4$ on the hyperbolic plane and constructing a Dirichlet polygon for the action.

\begin{thm}\label{PJ4 presentation}
The pure cactus group $PJ_4$ admits the following presentation: 
\[
      \left\langle 
      g_1, \cdots, g_{10} \left| 
      \begin{array}{l}
      g_{3} g_6^{-1} g_{7} g_{9}^{-1} g_{2}^{-1} 
      = g_{3} g_{8}^{-1} g_{4}^{-1}
      = g_{5} g_{9}^{-1} g_{4}^{-1} \\
      = g_{5} g_{1} g_{6}^{-1}
      = g_8 g_{10} g_7^{-1}
      = g_{10} g_{1}^{-1} g_{2} =e
      \end{array}
     \right.
     \right\rangle
\]
    
This presentation is transformed to the next one:
\[
     \left\langle g_2, g_4, g_8, g_9, g_{10} \mid g_{2}g_{9}g_{10}^{-1}g_{8}^{-1}g_{4}g_{9}g_{2}g_{10}g_{8}^{-1}g_{4}^{-1} =e
     \right\rangle
\]
\end{thm}

We remark that our generators $g_1, \dots, g_{10}$ are expressed as words in the standard generators of the cactus group $J_4$.
In fact, the generators $g_i$ described above have minimal length with respect to the standard generating set $\{ s_{i,j} \}$ of $J_4$; that is, their word lengths are either 4 or 5, and there are no elements of $PJ_4$ whose word length is 3 or less.
See Section~\ref{sec:proof} for details.

It can be verified that the presentation given in Theorem~\ref{PJ4 presentation} is equivalent to that of \cite[Appendix A]{BCL24}.

\begin{cor}[{\cite[Theorem~5.5]{BCL24}}]\label{BCL}
    The pure cactus group $PJ_4$ of degree four admits the following presentation. 
    \[
     \left\langle \alpha, \beta, \gamma, \delta, \epsilon \mid \alpha \gamma \epsilon \beta \epsilon \alpha^{-1} \delta^{-1} \beta \gamma \delta^{-1} =e\right\rangle
    \]
\end{cor}

As a corollary, we obtain a direct alternative proof of the fact that $PJ_4$ is isomorphic to the fundamental group of the connected sum of five real projective planes.

\begin{cor}\label{Pi_1_of_ConnSum}
    The pure cactus group $PJ_4$ of degree four admits the following presentation. 
\[
\left\langle 
\alpha_1, \alpha_2, \alpha_3, \alpha_4, \alpha_5 \mid 
\alpha_1^2 \alpha_2^2 \alpha_3^2 \alpha_4^2 \alpha_5^2 =e
\right\rangle.
\]
\end{cor}

The key to the proof of the above theorem is an action of $PJ_4$ on the Cayley complex of the subgroup $J_4^{[2,3]}$ of $J_4$, which is isometric to the hyperbolic plane $\mathbb{H}^2$ up to scaling.
We note that this action was originally introduced in \cite{genevois2022cactusgroupsviewpointgeometric}.

\begin{remark}
It is known that the moduli space $M_{0,n+1}(\mathbb{R})$ of real genus 0 curves over $\mathbb{R}$ with $n + 1$ marked points can be identified with the configuration space of distinct $n+1$ points on the circle.
See, for example, \cite{DJS03}.
The idea of our proof in this paper comes from the cell complex structure of such configuration spaces studied in \cite{AperyYoshida, MYos96, KNY99, MN00}.
Additionally, the relationship between the cactus group of degree 3 and such a configuration space is explored in \cite{HamaIchiharaPJ3}.
\end{remark}

\section{Cactus group}\label{subsec21}
In this section, we give definitions of the cactus group and the pure cactus group algebraically. 

For any integer $n \geq 2$, we define the \textit{cactus group of degree $n$}, denoted by $J_n$, as the group given by the presentation with generators $s_{p,q}$ for $1 \leq p < q \leq n$, subject to the following relations.
\begin{itemize}
   \item $s_{p,q}^2 = e$ for every $1 \leq p < q \leq n$,
   \item $s_{p,q}s_{m,r} = s_{m,r}s_{p,q}$ for all $1 \leq p < q \leq n$ and $1 \leq m < r \leq n$ satisfying $[p, q] \cap [m, r] = \emptyset$,
   \item $s_{p,q}s_{m,r} = s_{p+q-r,p+q-m}s_{p,q}$ for all $1 \leq p < q \leq n$ and $1 \leq m < r \leq n$ satisfying $[m, r] \subset [p, q]$.
\end{itemize}
Here, $e$ denotes the identity element, and $[p, q]$ denotes the set ${ p, p+1, \dots, q-1, q }$ of integers for positive integers $p$ and $q$ with $p < q$.
Throughout the sequel, we write generators such as $s_{23}$ (without a comma) instead of $s_{2,3}$ for brevity.

For the degree four case, the cuctus group $J_4$ has the following presentation. 
\[
\left\langle 
s_{12}, s_{23}, s_{34}, s_{13}, s_{24}, s_{14} 
\left| 
\begin{array}{l}
s_{12}^2= s_{23}^2= s_{34}^2= s_{13}^2= s_{24}^2= s_{14}^2=e,\\ 
s_{12} s_{34} = s_{34} s_{12}, 
s_{12} s_{13} = s_{13} s_{23}, 
s_{23} s_{24} = s_{24} s_{34}, \\
s_{12} s_{14} = s_{14} s_{34}, 
s_{23} s_{14} = s_{14} s_{23}, 
s_{13} s_{14} = s_{14} s_{24} 
\end{array}
\right.
\right\rangle.
\]

Similar to the braid group, elements of the cactus group can be represented by diagrams of vertical strands on the plane. 
Some such diagrams for $J_4$ are shown in Figure~\ref{Fig:element}.

\begin{figure}[htb]
    \begin{picture}(350,65)(0,0)
    \put(-20,0){
    \begin{tikzpicture}

\draw[ultra thick] (0,0) .. controls +(0,1) and +(0,-1) ..
(.5,2) ;
\draw[ultra thick] (.5,0) .. controls +(0,1) and +(0,-1) ..
(0,2) ;
\draw[ultra thick] (1,0) .. controls +(0,0) .. (1,2);
\draw[ultra thick] (1.5,0) .. controls +(0,0) .. (1.5,2);

\draw[ultra thick] (3,0) .. controls +(0,0) .. (3,2);
\draw[ultra thick] (3.5,0) .. controls +(0,1) and +(0,-1) ..
(4,2) ;
\draw[ultra thick] (4,0) .. controls +(0,1) and +(0,-1) ..
(3.5,2) ;
\draw[ultra thick] (4.5,0) .. controls +(0,0) .. (4.5,2);

\draw[ultra thick] (6,0) .. controls +(0,0) .. (6,2);
\draw[ultra thick] (6.5,0) .. controls +(0,0) .. (6.5,2);
\draw[ultra thick] (7,0) .. controls +(0,1) and +(0,-1) ..
(7.5,2) ;
\draw[ultra thick] (7.5,0) .. controls +(0,1) and +(0,-1) ..
(7,2) ;

\draw[ultra thick] (9,0) .. controls +(0,1) and +(0,-1) ..
(10,2) ;
\draw[ultra thick] (9.5,0) .. controls +(0,0) .. (9.5,2);
\draw[ultra thick] (10,0) .. controls +(0,1) and +(0,-1) ..
(9,2) ;
\draw[ultra thick] (10.5,0) .. controls +(0,0) .. (10.5,2);

\draw[ultra thick] (12,0) .. controls +(0,1) and +(0,-1) ..
(13.5,2) ;
\draw[ultra thick] (12.5,0) .. controls +(0,1) and +(0,-1) ..
(13,2) ;
\draw[ultra thick] (13,0) .. controls +(0,1) and +(0,-1) ..
(12.5,2) ;
\draw[ultra thick] (13.5,0) .. controls +(0,1) and +(0,-1) ..
(12,2) ;
\end{tikzpicture}
}
    \put(340,-10){$s_{14}$}
    \put(255,-10){$s_{13}$}
    \put(0,-10){$s_{12}$}
    \put(85,-10){$s_{23}$}
    \put(170,-10){$s_{34}$}

\end{picture}

%\vspace{1.5cm}

   \begin{picture}(350,110)(0,0)
    \put(20,0){
\begin{tikzpicture}
\draw[ultra thick] (0,0) .. controls +(0,1) and +(0,-1) ..
(1,1.5) .. controls +(0,1) and +(0,-1) .. (0,3);
\draw[ultra thick] (1,0) .. controls +(0,1) and +(0,-1) ..
(0,1.5) .. controls +(0,1) and +(0,-1) .. (1,3);
\draw[ultra thick] (.5,0) .. controls +(0,0) .. (.5,3);
\draw[ultra thick] (1.4,0) .. controls +(0,0) .. (1.4,3);

\draw[ultra thick] (2.5,0) .. controls +(0,0) .. (2.5,3);
\draw[ultra thick] (3,0) .. controls +(0,0) .. (3,3);
\draw[ultra thick] (3.5,0) .. controls +(0,0) .. (3.5,3);
\draw[ultra thick] (4,0) .. controls +(0,0) .. (4,3);

\draw[ultra thick] (6,0) .. controls +(0,1) and +(0,-1) ..
(6.5,1.5) .. controls +(0,1) and +(0,-1) .. (7,3);
\draw[ultra thick] (6.5,0) .. controls +(0,1) and +(0,-1) ..
(6,1.5) .. controls +(0,1) and +(0,-1) .. (7.5,3);
\draw[ultra thick] (7,0) .. controls +(0,0) ..
(7,1.5) .. controls +(0,1) and +(0,-1) .. (6.5,3);
\draw[ultra thick] (7.5,0) .. controls +(0,0) ..
(7.5,1.5) .. controls +(0,1) and +(0,-1) .. (6,3);

\draw[ultra thick] (9,0) .. controls +(0,1) and +(0,-1) ..
(10.5,1.5) .. controls +(0,1) and +(0,-1) .. (10,3);
\draw[ultra thick] (9.5,0) .. controls +(0,1) and +(0,-1) ..
(10,1.5) .. controls +(0,1) and +(0,-1) .. (10.5,3);
\draw[ultra thick] (10,0) .. controls +(0,1) and +(0,-1) ..
(9.5,1.5) .. controls +(0,0) .. (9.5,3);
\draw[ultra thick] (10.5,0) .. controls +(0,1) and +(0,-1) ..
(9,1.5) .. controls +(0,0) .. (9,3);

\end{tikzpicture}
}
    \put(77,40){\large $=$}
    \put(55,-15){$ ( s_{13} )^2 = e$}
    
    \put(255,40){\large $=$}
    \put(240, -15){$s_{14}s_{12} = s_{34} s_{14}$}
    % \put(9.2,1.2){$s_{12}s_{13}$}
\end{picture}

\vspace{.5cm}

\caption{Diagrams for some elements of $J_4$}\label{Fig:element}
\end{figure}

Due to this diagrammatic expression, we see that the cactus group $J_n$ admits a natural projection $\pi : J_n \to S_n$ onto the symmetric group $S_n$ of degree $n$, and its kernel is called the \textit{pure cactus group} of degree $n$, denoted by $PJ_n$. 
See \cite[Subsection 3.1]{HENRIQUES-KAMNITZER} or \cite[Section 1]{genevois2022cactusgroupsviewpointgeometric} for more details.

We also introduce a certain subgroup that will be used in later sections.
For each integer $n \geq 2$ and subset $S \subset [2, n]$, let $J_n^S$ be the subgroup of $J_n$ generated by the elements $s_{p,q}$ for $1 \leq p < q \leq n$ with $q - p + 1 \in S$, and defined by the following relations.
\begin{itemize}
   \item $s_{p,q}^2 = e$ for every $1 \leq p < q \leq n$ satisfying $q-p+1\in S$,
   \item $s_{p,q}s_{m,r} = s_{m,r}s_{p,q}$ for every $1 \leq p < q \leq n$ and $1 \leq m < r \leq n$ satisfying $[p, q] \cap [m, r] = \emptyset$ and $q-p+1\in S$,
   \item $s_{p,q}s_{m,r} = s_{p+q-r,p+q-m}s_{p,q}$ for every $1 \leq p < q \leq n$ and $1 \leq m < r \leq n$ satisfying $[m, r] \subset [p, q]$ and $q-p+1\in S$.
\end{itemize}
See \cite[Section 5]{genevois2022cactusgroupsviewpointgeometric} for more details. 

In later sections, we mainly consider the group $J_4^{[2,3]}$, which we denote by $J_4'$ in the rest of paper. 
It is a subgroup of $J_4$ generated by $s_{12}, s_{23}, s_{34}, s_{13}, s_{24} $ (excluding $s_{14}$), and it has the following presentation. 
\[
\left\langle 
s_{12}, s_{23}, s_{34}, s_{13}, s_{24} 
\left| 
\begin{array}{l}
s_{12}^2= s_{23}^2= s_{34}^2= s_{13}^2= s_{24}^2=e,\\ 
s_{12} s_{34} = s_{34} s_{12}, 
s_{12} s_{13} = s_{13} s_{23}, 
s_{23} s_{24} = s_{24} s_{34}
\end{array}
\right.
\right\rangle.
\]

We denote by $\mathcal{G}$ the Cayley graph of the group $J_4'$, and by $\mathcal{C}$ the Cayley complex of $J_4'$ with respect to the presentation above.
See, for instance, \cite{AT} for the definitions of these objects.
    
\section{Preliminary}

In this section, we prepare some properties of group actions for the proof of Theorem~\ref{PJ4 presentation}.

The following two propositions are given in \cite{Ratcliffe}.
See \cite{Ratcliffe} for the definitions of the terms used.

\begin{prop}[{\cite[Theorem~6.6.13]{Ratcliffe}}]\label{Dirichlet domain to fundamental domain}
    Let $D(a)$ be the Dirichlet domain with center $a$ for a discontinuous group $G$ of isometries of a metric space $X$ such that 
    \begin{enumerate}
        \item $X$ is geodesically connected, 
        \item $X$ is geodesically complete, 
        \item $X$ is finitely compact. 
    \end{enumerate}
    Then $D(a)$ is a locally finite fundamental domain for $G$. 
\end{prop}

\begin{prop}[{\cite[Exercise~6.6, 2.]{Ratcliffe}}]\label{D generates G}
    Let G be a group of isometries of a connected metric space $X$ with locally finite fundamental domain $R$. Then $G$ is generated by $\displaystyle\{g\in G \mid \overline{R} \cap g\overline{R} \neq \emptyset \}$. 
\end{prop}

Next, we recall the following well-known theorem, originally from \cite{Poincare}.
Our terminology follows that of \cite{Maskit}.

\begin{prop}[{Poincar\'{e}'s polygon theorem, c.f. \cite[p223]{Maskit}}]
\label{poincare's thm}
    Let $D$ be a Poincar\'{e} polygon. Let $G$ be the group generated by the identifying generators. Then $G$ is discontinuous, $D$ is a fundamental polygon for $G$, and the cycle relations form the complete set of relations for $G$. 
\end{prop}

In the following, we explain the terms used in the theorem above.
% In the following, we give explanations for the terms used in the theorem above. 

Let $D$ be a polygon in the hyperbolic plane $\mathbb{H}^2$.
An \textit{identification} on $D$ is a map that assigns, to each side $s$, a side $s'$ and an isometry $A(s, s')$ of $\mathbb{H}^2$ such that 
\begin{enumerate}
\item[(i)] $A(s, s')$ maps $s$ onto $s'$,
\item[(ii)] $(s')' = s$ and $A(s', s) = (A(s, s'))^{-1}$,
\item[(iii)] if $s = s'$, then $A(s, s')$ is the identity on $s$, and
\item[(iv)] for each side $s$, there exists a neighborhood $V$ of $s$ such that $A(s, s')(V \cap D) \cap D = \emptyset$.
\end{enumerate}
% Let $D$ be a polygon on the hyperbolic plane $\mathbb{H}^2$. 
% An \textit{identification} on $D$ is a map which assigns, to each side $s$, a side $s'$ and an isometry $A(s, s')$ of $\mathbb{H}^2$ so that
%     \begin{enumerate}
%         \item[(i)] $A(s,s')$ maps $s$ onto $s'$, 
%         \item[(ii)] $(s')' = s$ and $A(s',s) = ( A(s,s') )^{-1}$, 
%         \item[(iii)] if $s=s'$, then $A(s,s')$ is the identity on $s$, and 
%         \item[(iv)] for each side $s$, there exists a neighborhood $V$ of $s$ so that $A(s,s')(V \cap D) \cap D = \emptyset$. 
%     \end{enumerate}

We call the isometries $A(s, s')$ which identify the sides of the polygon $D$ \textit{identifying generators}, or \textit{generators}. 

Let $D$ be a polygon with identification. 
The \textit{identified polygon} $D^{*}$ is the quotient obtained by the surjection $p:\overline{D} \longrightarrow D^{*}$ such that $p (x) = p (x')$ if there exists a generator $A(s,s')$ with $A(s,s')(x) = x'$. 
For two elements $x$, $x'$ in $D^{*}$, we set a metric on $D^{*}$ as 
    \[
        d_{D^*}^*(x, x') := \inf \sum_{i=1}^n d_{\mathbb{H}^2}(z_i, z_i')
    \]
where the infimum is taken over all $n$ and over all $2n$-tuples of points of $\overline{D}$ such that $p(z_1) = x$, $p(z_i') = p(z_{i+1})$, and $p(z_n') = x'$. 

Let $D$ be a polygon with identification. 
The polygon $D$ is \textit{complete} if 
\begin{enumerate}
\item[(v)] for each $x \in D^{*}$, $p^{-1}(x)$ is a finite set, 
\item[(vi)] $D^{*}$ is complete in the metric $d^{*}_{\mathbb{H}^2}$. 
\end{enumerate}

Let $D$ be a complete polygon with identification. 
Let $z_1$ be a vertex of $D$. 
There are precisely two sides of $D$ meeting at $z_1$; we choose one of these and call it $s_1$.
There is then a corresponding side $s_1'$ and a generator $A_1 = A(s_1, s_1')$.
Set $z_2 = A_1(z_1)$ and observe that there is a unique other side $s_2$ having $z_2$ as an endpoint.
There is a corresponding side $s_2'$ and a generator $A_2 = A(s_2, s_2')$.
Set $z_3 = A_2(z_2)$ and observe that there is a unique other side $s_3$ having $z_3$ as an endpoint.
In this way, we obtain a sequence $\{z_i\}$ of vertices, and a sequence $\{ A_i \}$ of generators.
Then there exists a periodic sequence of generators $\{A_1, \dots, A_n\}$ such that
\[ A_n( \cdots A_2(A_1(z_1) ) ) = z_1 \] 
called a \textit{cycle of generators}. 
For a cycle of generators $\{A_1, \dots, A_n \}$, there exists a sequence of vertices $\{z_1, \dots, z_n\}$ in $D$ such that 
\[ z_2 = A_1(z_1) , \ \dots , \ z_n = A_{n-1}(z_{n-1}), \text{ and }  A_n(z_n) = z_1 . \] 
This sequence is called a \textit{cycle of vertices}. 
Then, the polygon $D$ is said to satisfy the \textit{cycle condition} if for each cycle of vertices $\{z_1, \dots, z_n\}$, there exists an integer $\nu$ such that 
    \begin{enumerate}
        \item[(vii)] $\displaystyle \nu \sum_{i=1}^{n} \alpha(z_i) = 2 \pi$, 
    \end{enumerate} 
where $\alpha(z_i)$ is the interior angle at $z_i$ measured inside $D$. 

Let $D$ be a polygon with identifications.  
The polygon $D$ is called a \textit{Poincar\'{e} polygon} if it is complete and satisfies the cycle condition.

For a cycle of generators $\{A_1, \dots, A_n\}$, by Condition (vii), the $\nu$-fold composition
\[
  (A_n \circ \cdots \circ A_1)^\nu
\]
is the identity on $D$, i.e., we have the relation
\[
  (A_n \circ \cdots \circ A_1)^\nu = e
\]
for that cycle.  
We call such relations the \textit{cycle relations}.

\section{Key propositions}\label{sec:prop}

In this section, we present two key propositions, Proposition~\ref{C4_32_isomH^2} and Proposition~\ref{prop42}, which are essential for the proof of Theorem~\ref{PJ4 presentation}.

\begin{prop}\label{C4_32_isomH^2}
The Caylay complex $\mathcal{C}$ of $J_4'$ is isometric to the hyperbolic plane $\mathbb{H}^2$ up to scaling. 
\end{prop}
\begin{proof}
As shown in \cite[Proposition 2.6]{genevois2022cactusgroupsviewpointgeometric}, every element of the cactus group has a unique normal form diagram. 
This enables us to enumerate all the elements of $J_4'$ with bounded word length, equivalently, to determine the local structure, vertices and edges, in $\mathcal{C}$ inside the bounded ball centered at $e$. 
Note that the vertices of $\mathcal{C}$ are identified with elements in $J_4'$.

The vertices in $\mathcal{C}$ whose word length is $1$ are simply the generators;
    \[
    \begin{array}{ccccc}
        s_{12} &
        s_{23} &
        s_{34} &
        s_{13} &
        s_{24}. 
    \end{array}
    \]
The following vertices have word length 2 with respect to the generators above.
    \[
    \begin{array}{ccccc}
        s_{12}s_{23}  &  s_{13}s_{23} & s_{13}s_{24} & s_{13}s_{34} & s_{13}s_{12} \\
        s_{23}s_{12} & s_{23}s_{34}  &  s_{24}s_{34} & s_{24}s_{12} & s_{24}s_{13} \\
        s_{24}s_{23} & s_{34}s_{23} & s_{34}s_{13}  &  s_{34}s_{12} & s_{12}s_{24}
    \end{array}
    \]

As shown in Figure~\ref{word_length_2}, there are 25 edges (1-cells) connecting the above vertices (0-cells) in $\mathcal{C}$, and the 2-cells (quadrangles) containing $e$ in $\mathcal{C}$ are
\[
\begin{array}{l}
\langle e, s_{12}, s_{12}s_{13}, s_{13} \rangle, 
\langle e, s_{13}, s_{13}s_{12}, s_{23} \rangle, 
\langle e, s_{23}, s_{23}s_{24}, s_{24} \rangle, \\
\langle e, s_{24}, s_{24}s_{23}, s_{34} \rangle, 
\langle e, s_{34}, s_{34}s_{12}, s_{12} \rangle.
\end{array}
\]

\begin{figure}[htb]
    \centering
\begin{overpic}[width=.6\textwidth]{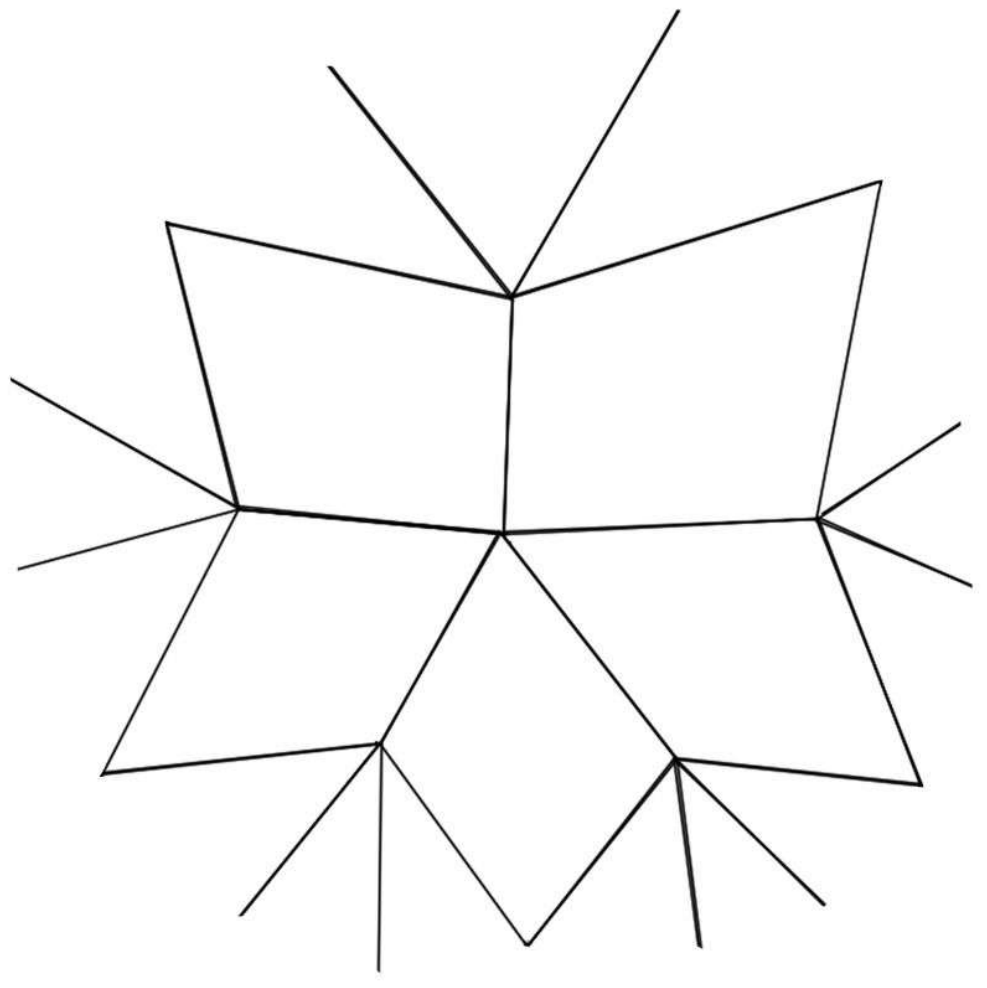}
\put(33,51){$e$} 

    \put(49,52){$s_{12}$}
    \put(31,61){$s_{13}$}
    \put(20,53){$s_{23}$}
    \put(24,39){$s_{24}$} 
    \put(44,38){$s_{34}$}

        \put(61,57){$s_{12}s_{23}$}
        \put(55,71){$s_{13}s_{23}$}
        \put(45,81){$s_{13}s_{24}$}
        \put(25,78){$s_{13}s_{34}$}
        \put(11,69){$s_{13}s_{12}$}
        \put(2,60){$s_{23}s_{12}$}
        \put(1,49){$s_{23}s_{34}$}
        \put(8,34){$s_{24}s_{34}$}
        \put(15,26){$s_{24}s_{12}$}
        \put(24,22){$s_{24}s_{13}$}
        \put(35,24){$s_{24}s_{23}$}
        \put(45,24){$s_{34}s_{23}$}
        \put(53,27){$s_{34}s_{13}$}
        \put(59,36){$s_{34}s_{12}$}
        \put(62,47){$s_{12}s_{24}$}       
\end{overpic}

\vspace{-2cm}

    \caption{Neighborhood of $e$ in the Cayley complex $\mathcal{C}$}
    \label{word_length_2}
\end{figure}

Therefore, $\mathcal{C}$ is a cell complex such that there are five quadrilateral $2$-cells surrounding each vertex.

We consider the map $f$ corresponding to each closed $2$-cell on $\mathcal{C}$ to the square (regular quadrilateral) in $\mathbb{H}^2$ whose internal angles are equal to $\frac{2 \pi}{5}$.
The image of the map $f$ gives a uniform tiling (regular tessellation) of $\mathbb{H}^2$ of type $\{4,5\}$.
See Figure~\ref{fig:45tesse}.

\begin{figure}[htb]
\includegraphics[width=.6\textwidth]{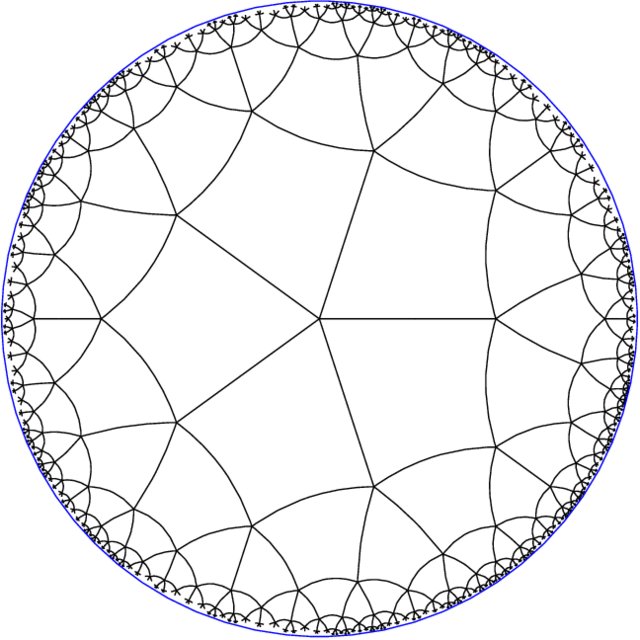}
\caption{The $\{4,5\}$-tesselation of the hyperbolic plane. (\cite[Figure 7]{BandeltChepoiEppstein})}
\label{fig:45tesse}
\end{figure}
 
We see that this map is an isometry up to scaling, satisfying $d(f(x), f(y)) = R \, d(x, y)$ with $R = \cosh^{-1}\left(\cot^2{\frac{\pi}{5}}\right)$ by an elementary calculation based on hyperbolic trigonometry.
\end{proof}

In other words, if we set a metric on $\mathcal{C}$ by assigning each edge a constant length $R$, then $\mathcal{C}$ is truly isometric to $\mathbb{H}^2$.  
Therefore, we can apply results for the hyperbolic plane directly to the Cayley complex $\mathcal{C}$.

Next, we consider an action of $PJ_4$ on the Cayley complex $\mathcal{C}$. 
The following proposition was essentially obtained in \cite{genevois2022cactusgroupsviewpointgeometric}. 

\begin{prop}[{\cite[Proof of Corollary~7.3]{genevois2022cactusgroupsviewpointgeometric}}]\label{prop42}
Let $\Gamma_{0}$ be the action of $PJ_4$ on $J_4'$ defined as follows. 
    \begin{align*}
        \Gamma_{0}: PJ_4 \times J_4'
        &\longrightarrow J_4'\\
        (g, h) &\longmapsto
        \begin{cases}
            gh &  gh \in J_4'\\
            ghs_{14} & gh \notin J_4'
        \end{cases}
    \end{align*}
Then $\Gamma_0$ induces the maps $\Gamma_1$ and $\Gamma$, which give rise to isometric actions of $PJ_4$ on $\mathcal{G}$ and $\mathcal{C}$, respectively.  
Note that $\left( \mathcal{C} \right)^{(0)} = J_4'$ is identified with $J_4'$. 
Moreover, the actions $\Gamma_1$ and $\Gamma$ on $\mathcal{G}$ and $\mathcal{C}$ are free and cocompact. 
\end{prop}

In the following, for simplicity, we identify $PJ_4$ with its image under the action $\Gamma$, and write $g \cdot x$ instead of $\Gamma(g, x)$ for $g \in PJ_4$ and $x \in \mathcal{C}$.  
In this setting, $PJ_4$ is regarded as a subgroup of $\mathrm{Isom}\;\mathcal{C}$.

\section{Proof of Theorem~\ref{PJ4 presentation}}\label{sec:proof}

In this section, we give our proof of Theorem~\ref{PJ4 presentation}.  
First, we will find a fundamental polygon $\widetilde{D}$ on $\mathcal{C}$, equivalently on $\mathbb{H}^2$, with respect to the action $\Gamma$ of $PJ_4$, and apply Poincar\'{e}'s theorem to obtain a presentation of $PJ_4$.

We first prepare the next lemma. 

\begin{lemma}\label{g_i are in PJ4}
Let $a_1, \cdots, a_{20}$ be the elements of $J_4'$ defined as follows. 
\[\begin{array}{lll}
        a_1 := s_{13}s_{24}s_{12}s_{34} \qquad&
        a_2 := s_{13}s_{24}s_{13}s_{24} \qquad&
        a_3 := s_{13}s_{34}s_{23}s_{12} \\
        a_4 := s_{13}s_{34}s_{13}s_{23} &
        a_5 := s_{23}s_{12}s_{23}s_{13} &
        a_6 := s_{23}s_{12}s_{24}s_{12} \\
        a_7 := s_{23}s_{34}s_{13}s_{34} &
        a_8 := s_{24}s_{34}s_{23}s_{34} &
        a_9 := s_{24}s_{12}s_{24}s_{23} \\
        a_{10} := s_{24}s_{12}s_{23}s_{34} &
        a_{11} := s_{24}s_{13}s_{24}s_{13} &
        a_{12} := s_{24}s_{23}s_{13}s_{34} \\
        a_{13} := s_{34}s_{23}s_{34}s_{24} &
        a_{14} := s_{34}s_{23}s_{12}s_{24} &
        a_{15} := s_{34}s_{13}s_{34}s_{23} \\
        a_{16} := s_{34}s_{13}s_{23}s_{24} &
        a_{17} := s_{34}s_{12}s_{24}s_{13} &
        a_{18} := s_{12}s_{24}s_{12}s_{23} \\
        a_{19} := s_{12}s_{23}s_{34}s_{13} &
        a_{20} := s_{12}s_{23}s_{12}s_{13} &
\end{array}\]

    Then the following elements $g_1, \dots , g_{10}$ and their inverses $g_1^{-1}, \dots , g_{10}^{-1}$ are contained in $PJ_4$.     
\[
\begin{array}{ll}
    g_{1} = a_{1}s_{14}    \qquad  &g_{1}^{-1} = a_{16}s_{14}\\
    g_{2} = a_{2}           &g_{2}^{-1} = a_{11}\\
    g_{3} = a_{3}s_{14}     &g_{3}^{-1} = a_{14}s_{14}\\
    g_{4} = a_{4}s_{14}     &g_{4}^{-1} = a_{9}s_{14}\\
    g_{5} = a_{5}           &g_{5}^{-1} = a_{20}\\ 
    g_{6} = a_{6}s_{14}     &g_{6}^{-1} = a_{15}s_{14}\\
    g_{7} = a_{7}s_{14}     &g_{7}^{-1} = a_{18}s_{14}\\
    g_{8} = a_{8}           &g_{8}^{-1} = a_{13}\\
    g_{9} = a_{10}s_{14}    &g_{9}^{-1} = a_{19}s_{14}\\
    g_{10} = a_{12}s_{14}   &g_{10}^{-1} = a_{17}s_{14}
\end{array}
\]

Furthermore, with respect to the action $\Gamma_1$ on $\mathcal{G}$ in Proposition~\ref{prop42}, they are all of the elements of $PJ_4$ whose translation lengths are at most 4 and actually equal to $4$. 
Precisely, if $ d( e, g \cdot e) \le 4 $
%    \[
%    \min\{d(x,gx) \mid x \in J_4'\} = 4
%    \]
holds for $g$ in $PJ_4$, then $g = g_i^{\pm1}$ and $ d( e, g_i \cdot e) = 4 $ for some $1 \leq i \leq 10$. 
\end{lemma}

\begin{proof}
We first confirm that $g_1, ..., g_{10}$ and their inverses $g_1^{-1}, \dots , g_{10}^{-1}$ are in $PJ_4$. 
The elements $a_1, \cdots, a_{8}, a_{10}, a_{12}$ correspond to the elements in $S_4$ via $\pi$ as follows. 
\begin{align*}
&\pi(a_{1}) = \pi(s_{13}s_{24}s_{12}s_{34})
=(13)(24)(12)(34)
=(14)(23)\\
&\pi(a_{2}) = \pi(s_{13}s_{24}s_{12}s_{34})
=(13)(24)(13)(24)
=e\\
&\pi(a_{3}) = \pi(s_{13}s_{34}s_{23}s_{12})
=(13)(34)(23)(12)
=(14)(23)\\
&\pi(a_{4}) = \pi(s_{13}s_{34}s_{13}s_{23})
=(13)(34)(13)(23)
=(14)(23)\\
&\pi(a_{5}) = \pi(s_{23}s_{12}s_{23}s_{13})
=(23)(12)(23)(13)
=e\\
&\pi(a_{6}) = \pi(s_{23}s_{12}s_{24}s_{12})
=(23)(12)(24)(12)
=(14)(23) \\
&\pi(a_{7}) = \pi(s_{23}s_{34}s_{13}s_{34})
=(23)(34)(13)(34)
=(14)(23) \\
&\pi(a_{8}) = \pi(s_{24}s_{34}s_{23}s_{34})
=(24)(34)(23)(34)
=e \\
&\pi(a_{10}) = \pi(s_{24}s_{12}s_{23}s_{34})
=(24)(12)(23)(34)
=(14)(23) \\
&\pi(a_{12}) = \pi(s_{24}s_{23}s_{13}s_{34})
=(24)(23)(13)(34)
=(14)(23) 
\end{align*} 
Thus, by definition, the elements $g_1, ..., g_{10}$ are in $PJ_4$. 

Next we show their inverses $g_1^{-1}, \dots , g_{10}^{-1}$ can be described as in the statement of the lemma. 
For $a_1 , \dots , a_{20}$, the following are obtained by straightforward calculations. 
\[
\begin{array}{llll}
            a_{1}^{-1} = a_{17} \qquad &
            a_{2}^{-1} = a_{11} \qquad &
            a_{3}^{-1} = a_{19} \qquad &
            a_{4}^{-1} = a_{4} \\
            a_{5}^{-1} = a_{20} \qquad &
            a_{6}^{-1} = a_{18} \qquad &
            a_{7}^{-1} = a_{15} \qquad &
            a_{8}^{-1} = a_{13} \\
            a_{9}^{-1} = a_{9} \qquad &
            a_{10}^{-1} = a_{14} \qquad &
            a_{12}^{-1} = a_{16} & 
\end{array}
\]
For example, we see; 
\[
\begin{array}{ll}
a_{1}a_{17} = s_{13}s_{24}s_{12}s_{34}s_{34}s_{12}s_{24}s_{13} =e ,\\
a_{17}a_{1} = s_{34}s_{12}s_{24}s_{13}s_{13}s_{24}s_{12}s_{34} =e. 
\end{array}
\]
We can also obtain 
\[
\begin{array}{ll}
    s_{14}a_{1} = a_{12}s_{14} \qquad &
    s_{14}a_{2} = a_{11}s_{14} \\
    s_{14}a_{3} = a_{10}s_{14} \qquad &
    s_{14}a_{4} = a_{9}s_{14} \\
    s_{14}a_{5} = a_{8}s_{14} \qquad &
    s_{14}a_{6} = a_{7}s_{14} \\
    s_{14}a_{7} = a_{6}s_{14} \qquad &
    s_{14}a_{8} = a_{5}s_{14} \\
    s_{14}a_{9} = a_{4}s_{14} \qquad &
    s_{14}a_{10} = a_{3}s_{14} \\
    s_{14}a_{11} = a_{2}s_{14} \qquad &
    s_{14}a_{12} = a_{1}s_{14} 
\end{array}
\]
by simple calculation in the same way as the following. 
\begin{align*}
    s_{14}a_{1} &=s_{14}s_{13}s_{24}s_{12}s_{34}
    =s_{24}s_{14}s_{24}s_{12}s_{34}\\  
    &=s_{24}s_{13}s_{14}s_{12}s_{34} 
    =s_{24}s_{13}s_{34}s_{14}s_{34}\\ 
    &=s_{24}s_{13}s_{34}s_{12}s_{14}
    =s_{24}s_{13}s_{12}s_{34}s_{14}\\ 
    &=s_{24}s_{23}s_{13}s_{34}s_{14}
    =a_{12}s_{14} 
\end{align*}
Thus, when we set $g_1, \dots, g_{20}$ as in the statement, we have 
\[
\begin{array}{ll}
    g_{1}^{-1} = a_{16}s_{14} \qquad 
    &g_{2}^{-1} = a_{11}\\
    g_{3}^{-1} = a_{14}s_{14}
    &g_{4}^{-1} = a_{9}s_{14}\\
    g_{5}^{-1} = a_{20}
    &g_{6}^{-1} = a_{15}s_{14}\\
    g_{7}^{-1} = a_{18}s_{14}
    &g_{8}^{-1} = a_{13}\\
    g_{9}^{-1} = a_{19}s_{14}
    &g_{10}^{-1} = a_{17}s_{14}. 
\end{array}
\]
by 
\begin{align*}
    g_{1}^{-1} 
    &= ( a_1 s_{14} )^{-1} 
    = ( s_{14} a_{12} )^{-1} \\
    &= a_{12}^{-1} s_{14} 
    = a_{16} s_{14} 
\end{align*}
and so on. 

Finally, with respect to the action $\Gamma_1$ on $\mathcal{G}$ in Proposition~\ref{prop42}, we prove that $g_1, \dots, g_{10}$ and their inverses are all the elements of $PJ_4$ whose translation lengths equal $4$, i.e., if $d(e, g \cdot e) \le 4$ holds for some $g \in PJ_4$, then $g = g_i^{\pm 1}$ and $d(e, g_i \cdot e) = 4$ for some $1 \leq i \leq 10$.

Suppose that an element $g$ of $PJ_4$ has translation length at most $4$, i.e., $d(e, g \cdot e) \le 4$ holds.  
This means that the element $g \cdot e$ in $J_4'$ has word length at most $4$.  

All the elements of $J_4'$ whose word lengths are at most $4$ can be listed by brute-force enumeration. 
See Tables~\ref{tab:length3} and \ref{tab:length4}.  
In fact, there are $40$ elements of length $3$ and $105$ elements of length $4$.  

If $a := g \cdot e$ in $J_4'$, then, by the next claim, it lies in the image of $PJ_4$ under the map $\mathfrak{p}$, and so $\pi(a) \in \{e, (14)(23)\}$.

\begin{claim}\label{(a_i)s_14 is in PJ_4}
Let $\pi$ denote the natural surjection from $J_4$ to $S_4$ and $\mathfrak{p} : PJ_4 \to J_4'$ the map defined by 
\begin{align*}
    g \longmapsto        
        \begin{cases}    
            g &  \text{ if } \ g \in J_4'\\
            gs_{14} &  \text{ if } \ g \notin J_4'.
        \end{cases}
\end{align*}        
Then the following holds:
\[
    \mathfrak{p}(PJ_4) = \left( \pi \vert_{J_4'} \right)^{-1}(\{e, (14)(23)\}).
\]
\end{claim}

\begin{proof} 
For any element $x$ in $\left( \pi \vert _{J_4'} \right )^{-1}(\{e, (14)(23)\})$, we have $\pi(x) = e$ or $(14)(23)$.  
Since $x$ lies in $J_4'$, it does not contain $s_{14}$ as a word.  
In the former case, $x$ is in $PJ_4$ by definition, and thus $\mathfrak{p}(x) = x$.  
In the latter case, the following equalities hold.
    \[
    \begin{array}{ll}
         \pi(x) &= (14)(23) \\
         \pi(x) \pi(s_{14}) &= (14)(23)(14)(23) \\
         \pi(x s_{14}) &= e 
    \end{array}
    \]
Thus, $x s_{14}$ is in $PJ_4$ and $x s_{14}$ is not in $J_4'$. 
It follows that $\mathfrak{p}(x s_{14}) = x s_{14}^2 = x$.  

Conversely, for each element $x$ in $\mathfrak{p}(PJ_4)$, we have $x \in J_4'$.  
By the construction of $\mathfrak{p}$, there exists $y \in PJ_4$ such that $x = y$ or $x = y s_{14}$.  
It follows that  
\[
    \pi \vert _{J_4'}(x) = e \quad \text{or} \quad \pi \vert _{J_4'}(x) = (14)(23).
\]  
Therefore, $x \in \left( \pi \vert _{J_4'} \right )^{-1}(\{e, (14)(23)\})$.
\end{proof}

Among the elements in Tables~\ref{tab:length3} and \ref{tab:length4}, by directly checking whether the images of the elements under $\pi$ are $e$ or $(14)(23)$, we obtain only the elements $a_1, \cdots, a_{20}$ as given above.  
(Remark that $a_{16} = s_{34} s_{13} s_{23} s_{24}$ appears as $s_{34} s_{13} s_{24} s_{34}$ in Table~\ref{tab:length4}.)  

Consequently, all elements of $PJ_4$ whose translation lengths are at most $4$ are shown to be $g_1, \dots, g_{10}$ and their inverse elements $g_1^{-1}, \dots, g_{10}^{-1}$.
\end{proof}

\begin{table}[htb]
    \[
    \begin{array}{cccccc}
s_{13}s_{23}s_{34} &
s_{13}s_{23}s_{24} &
s_{13}s_{24}s_{12} &
s_{13}s_{24}s_{13} &
s_{13}s_{24}s_{23} &
s_{13}s_{34}s_{23} \\
s_{13}s_{34}s_{13} &
s_{13}s_{34}s_{12} &
s_{23}s_{13}s_{24} &
s_{23}s_{13}s_{23} &
s_{23}s_{12}s_{23} &
s_{23}s_{12}s_{24} \\
s_{23}s_{12}s_{34} &
s_{23}s_{34}s_{13} &
s_{23}s_{34}s_{23} &
s_{23}s_{34}s_{24} &
s_{24}s_{34}s_{13} &
s_{24}s_{34}s_{12} \\
s_{24}s_{12}s_{24} &
s_{24}s_{12}s_{23} &
s_{24}s_{12}s_{13} &
s_{24}s_{13}s_{24} &
s_{24}s_{13}s_{34} &
s_{24}s_{13}s_{12} \\
s_{34}s_{24}s_{12} &
s_{34}s_{24}s_{34} &
s_{34}s_{23}s_{34} &
s_{34}s_{23}s_{12} &
s_{34}s_{23}s_{13} &
s_{34}s_{13}s_{34} \\
s_{34}s_{13}s_{24} &
s_{34}s_{13}s_{23} &
s_{12}s_{34}s_{23} &
s_{12}s_{34}s_{24} &
s_{12}s_{24}s_{13} &
s_{12}s_{24}s_{12} \\
s_{12}s_{24}s_{34} &
s_{12}s_{23}s_{34} &
s_{12}s_{23}s_{12} &
s_{12}s_{23}s_{13} & & 
    \end{array}\]
    \caption{Elements of $J_4'$ of length 3}
    \label{tab:length3}
\end{table}

\begin{table}[htb]
    \[
    \begin{array}{cccc}
s_{13}s_{23}s_{34}s_{12} &
s_{13}s_{23}s_{34}s_{13} &
s_{13}s_{23}s_{34}s_{23} &
s_{13}s_{23}s_{34}s_{24} \\
s_{13}s_{24}s_{34}s_{13} &
s_{13}s_{24}s_{34}s_{12} &
s_{13}s_{24}s_{12}s_{24} &
s_{13}s_{24}s_{12}s_{13}\\
s_{13}s_{24}s_{12}s_{13} &
s_{13}s_{24}s_{13}s_{24} &
s_{13}s_{24}s_{13}s_{34} &
s_{13}s_{24}s_{13}s_{12}\\
s_{13}s_{34}s_{24}s_{12} &
s_{13}s_{34}s_{24}s_{34} &
s_{13}s_{34}s_{23}s_{34} &
s_{13}s_{34}s_{23}s_{12}\\
s_{13}s_{34}s_{23}s_{13} &
s_{13}s_{34}s_{13}s_{34} &
s_{13}s_{34}s_{13}s_{24} &
s_{13}s_{34}s_{13}s_{23}\\
s_{13}s_{12}s_{34}s_{23} &
s_{13}s_{12}s_{34}s_{24} &
s_{23}s_{13}s_{24}s_{13} &
s_{23}s_{13}s_{24}s_{12}\\
s_{23}s_{13}s_{24}s_{34} &
s_{23}s_{12}s_{13}s_{34} &
s_{23}s_{12}s_{13}s_{12} &
s_{23}s_{12}s_{23}s_{12}\\
s_{23}s_{12}s_{23}s_{34} &
s_{23}s_{12}s_{23}s_{24} &
s_{23}s_{12}s_{24}s_{12} &
s_{23}s_{12}s_{24}s_{13}\\
s_{23}s_{12}s_{24}s_{23} &
s_{23}s_{34}s_{12}s_{23} &
s_{23}s_{34}s_{12}s_{13} &
s_{23}s_{34}s_{13}s_{24}\\
s_{23}s_{34}s_{13}s_{34} &
s_{23}s_{34}s_{13}s_{12} &
s_{23}s_{34}s_{23}s_{12} &
s_{23}s_{34}s_{23}s_{34}\\
s_{23}s_{34}s_{23}s_{24} &
s_{23}s_{24}s_{23}s_{12} &
s_{23}s_{24}s_{23}s_{13} &
s_{24}s_{34}s_{13}s_{34}\\
s_{24}s_{34}s_{13}s_{24} &
s_{24}s_{34}s_{13}s_{23} &
s_{24}s_{12}s_{34}s_{23} &
s_{24}s_{12}s_{34}s_{24}\\
s_{24}s_{12}s_{24}s_{13} &
s_{24}s_{12}s_{24}s_{12} &
s_{24}s_{12}s_{24}s_{34} &
s_{24}s_{12}s_{23}s_{34}\\
s_{24}s_{12}s_{23}s_{12} &
s_{24}s_{12}s_{23}s_{13} &
s_{24}s_{13}s_{23}s_{34} &
s_{24}s_{13}s_{23}s_{24}\\
s_{24}s_{13}s_{24}s_{12} &
s_{24}s_{13}s_{24}s_{13} &
s_{24}s_{13}s_{24}s_{23} &
s_{24}s_{13}s_{34}s_{23}\\
s_{24}s_{13}s_{34}s_{13} &
s_{24}s_{13}s_{34}s_{12} &
s_{24}s_{23}s_{13}s_{24} &
s_{24}s_{23}s_{13}s_{23}\\
s_{34}s_{24}s_{12}s_{23} &
s_{34}s_{24}s_{12}s_{24} &
s_{34}s_{24}s_{12}s_{34} &
s_{34}s_{23}s_{24}s_{13}\\
s_{34}s_{23}s_{24}s_{23} &
s_{34}s_{23}s_{34}s_{23} &
s_{34}s_{23}s_{34}s_{13} &
s_{34}s_{23}s_{34}s_{12}\\
s_{34}s_{23}s_{12}s_{24} &
s_{34}s_{23}s_{12}s_{23} &
s_{34}s_{23}s_{12}s_{13} &
s_{34}s_{13}s_{12}s_{24}\\
s_{34}s_{13}s_{12}s_{34} &
s_{34}s_{13}s_{34}s_{13} &
s_{34}s_{13}s_{34}s_{23} &
s_{34}s_{13}s_{34}s_{24}\\
s_{34}s_{13}s_{24}s_{13} &
s_{34}s_{13}s_{24}s_{12} &
s_{34}s_{13}s_{24}s_{34} &
s_{34}s_{12}s_{13}s_{34}\\
s_{34}s_{12}s_{13}s_{12} &
s_{12}s_{34}s_{23}s_{12} &
s_{12}s_{34}s_{23}s_{34} &
s_{12}s_{34}s_{23}s_{24}\\
s_{12}s_{24}s_{23}s_{12} &
s_{12}s_{24}s_{23}s_{13} &
s_{12}s_{24}s_{13}s_{34} &
s_{12}s_{24}s_{13}s_{24}\\
s_{12}s_{24}s_{13}s_{23} &
s_{12}s_{24}s_{12}s_{23} &
s_{12}s_{24}s_{12}s_{24} &
s_{12}s_{24}s_{12}s_{34}\\
s_{12}s_{23}s_{24}s_{13} &
s_{12}s_{23}s_{24}s_{23} &
s_{12}s_{23}s_{34}s_{23} &
s_{12}s_{23}s_{34}s_{13}\\
s_{12}s_{23}s_{34}s_{12} &
s_{12}s_{23}s_{12}s_{24} &
s_{12}s_{23}s_{12}s_{23} &
s_{12}s_{23}s_{12}s_{13}\\
s_{12}s_{13}s_{12}s_{24} &
%s_{12}s_{13}s_{12}s_{34} 
& & 
    \end{array}
    \]
    \caption{Elements of $J_4'$ of length 4}
    \label{tab:length4}
\end{table}

\begin{figure}[htb]
    \centering
    %\hspace{-2.4cm}
    %\vspace{-1cm}
\begin{overpic}[width=\textwidth]{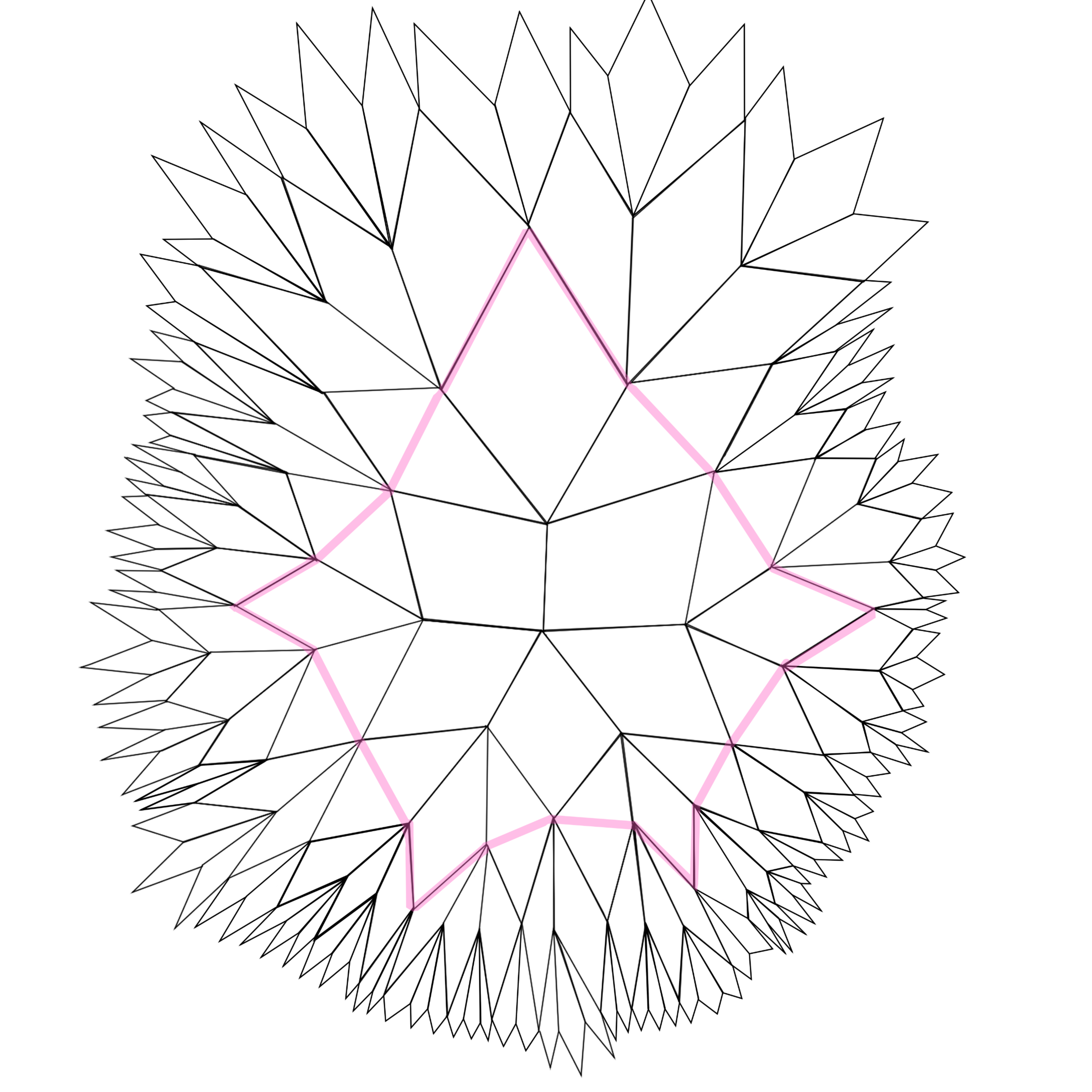}

\put(50,41){$e$} 

    \put(61,42){$s_{12}$}
    \put(50,50){$s_{13}$}
    \put(38,42){$s_{23}$}
    \put(45,32.5){$s_{24}$} 
    \put(57,32){$s_{34}$}

        \put(68,47){$s_{12}s_{23}$}
        \put(63.5,57){$s_{13}s_{23}$}
        \put(55.5,65){$s_{13}s_{24}$}
        \put(38,65){$s_{13}s_{34}$}
        \put(33,55.2){$s_{13}s_{12}$}
        \put(27,49){$s_{23}s_{12}$}
        \put(25,40.5){$s_{23}s_{34}$}
        \put(31.5,33){$s_{24}s_{34}$}
        \put(35,25){$s_{24}s_{12}$}
        \put(42,24){$s_{24}s_{13}$}
        \put(48.3,26){$s_{24}s_{23}$}
        \put(56,25){$s_{34}s_{23}$}
        \put(63,26){$s_{34}s_{13}$}
        \put(65,32){$s_{34}s_{12}$}
        \put(70,39){$s_{12}s_{24}$}
        
            \put(45,79){\small{$s_{13}s_{24}s_{23}$}}
            \put(20,45){\small{$s_{23}s_{34}s_{12}$}}
            \put(35,17){\small{$s_{24}s_{13}s_{23}$}}
            \put(62,19){\small{$s_{34}s_{13}s_{12}$}}
            \put(76,45){\small{$s_{12}s_{24}s_{34}$}}
\end{overpic}
    \caption{}
    \label{derichlet}
\end{figure}

Let $\widetilde{D}$ be the polygon on $\mathcal{C}$ defined as follows, which is visualized in Figure~\ref{derichlet}. 
\begin{align*}
    \displaystyle \widetilde{D}:=\bigcap_{i=1}^{20} \left\{x \in \mathcal{C} \mid \ d(e,x) \leq d(x,a_i) \right\} 
\end{align*}

From Figure~\ref{derichlet}, the following can be verified.

\begin{lemma}\label{v_1...v_5 in tilD}
With respect to the action $\Gamma_1$ on $\mathcal{G}$, the word lengths of the vertices in $\mathrm{int}\,\widetilde{D}$ are at most $3$.
The vertices on $\partial\widetilde{D}$ whose word length is $3$ are only 
$s_{13}s_{24}s_{23}$, 
$s_{23}s_{34}s_{12}$, 
$s_{24}s_{13}s_{23}$, 
$s_{34}s_{13}s_{12}$, 
$s_{12}s_{24}s_{34}$. \qed
\end{lemma} 
We denote the five vertices above by $v_1$, $\dots$, $v_5$, respectively.
We remark that there are ten sides of $\widetilde{D}$ which are not edges in the Cayley complex $\mathcal{C}$.  
These sides are the diagonals of some 2-cells (quadrangles) in $\mathcal{C}$.  
The interior angles at $v_1$, $\dots$, $v_5$ are $\frac{2\pi}{5}$.  
The interior angles at the vertices $s_{13}s_{23}$, $s_{13}s_{12}$, $s_{24}s_{34}$, $s_{24}s_{23}$, $s_{34}s_{12}$ on $\partial\widetilde{D}$ are $\frac{4\pi}{5}$.  
The interior angles at $s_{12}s_{23}$, $s_{13}s_{24}$, $s_{13}s_{34}$, $s_{23}s_{12}$, $s_{23}s_{34}$, $s_{24}s_{12}$, $s_{24}s_{13}$, $s_{34}s_{23}$, $s_{34}s_{13}$, $s_{12}s_{24}$ on $\partial\widetilde{D}$ are $\frac{3\pi}{5}$.

%See Figure~\ref{derichlet}. 
%\end{remark}

%\input{dirichlet_domain}

\begin{lemma}\label{D tilde is Dirichlet domain}
The subset $\widetilde{D}_1:= \mathrm{int}(\widetilde{D} \cap \mathcal{G})$ of $\mathcal{G}$ is the Dirichlet domain with center $e$ in $\mathcal{G}$ with respect to the action $\Gamma_1$ of $PJ_4$.% as a subgroup of $\mathrm{Isom}\;\mathcal{C}$. 
\end{lemma}

%In the context of a Cayley graph and its relation to group actions, a fundamental region (or fundamental domain) is a subset of the graph that, when acted upon by the group, covers the entire graph without overlaps (or with overlaps only on the boundary). 

\begin{proof}
Since $PJ_4$ is nontrivial, by definition, we consider 
\[
    D(e) = \bigcap_{g \neq e} H_g (e). 
\]
Here, $H_g (e)$ denotes the half-space defined by $\{x \in \mathcal{G} \mid d(x,e) < d(x,g \cdot e) \}$ for $g$ in $PJ_4$. 

We prove that $\widetilde{D}_1$ coincides with $D(e)$.  
By the definition of the Dirichlet domain, $\widetilde{D}_1$ contains $D(e)$.  
Thus, it suffices to show the converse, i.e., that $H_g(e)$ contains $\widetilde{D}_1$ for any nontrivial element $g$ in $PJ_4$.

Let $g$ be a nontrivial element of $PJ_4$.  
By Lemma~\ref{g_i are in PJ4}, if $g$ satisfies $d(e, g \cdot e) \le 4$, then $d(e, g \cdot e) = 4$ and $g = g_i^{\pm 1}$ for some $i$.  
Therefore, $H_g(e)$ is equal to $H_{g_i^{\pm 1}}(e)$ for some $i$ when $d(e, g \cdot e) \le 4$.

If $d(e, g \cdot e) \ge 6$, we obtain the following:
\[
    \{ x \mid d(e, x) < 3 \} \subset H_g(e).
\]
Since the word lengths of the vertices in $\widetilde{D}_1$ are less than $3$, 
we have $d(e, x) < 3$ for any $x \in \widetilde{D}_1$.  
Thus, 
\[
    \bigcap \{ H_g(e) \mid d(e, g \cdot e) \ge 6 \} \supset \widetilde{D}_1.
\]

It remains to show that $d(e, g \cdot e) \ne 5$ on $\mathcal{G}$ for any $g \in PJ_4$.
Recall that the images of the standard generators 
$s_{12}$, $s_{23}$, $s_{34}$, $s_{13}$, $s_{24}$ 
of $J_4'$ under $\pi$ are the transpositions 
$(12)$, $(23)$, $(34)$, $(13)$, $(24)$ in $S_4$, 
and also that $\pi(s_{14}) = (14)(23)$.  
See Section~\ref{subsec21}.

Let $g$ be an element in $PJ_4$. 
This means that $\pi(g) = e$ in $S_4$. 

If $g$ is contained in $J_4'$, then $g \cdot e = g$ is in $J_4'$. 
Thus $d(e, g \cdot e) = d(e, g ) = 5$ implies that $g$ has word length 5 with respect to the generators 
$s_{12}$, $s_{23}$, $s_{34}$, $s_{13}$, $s_{24}$ in $J_4'$. 
This would imply that $\pi(g)$ is an odd permutation, contradicting $\pi(g) = e$. 

If $g$ is not contained in $J_4'$, then $g \cdot e = g s_{14}$ is in $J_4'$ by definition of the action $\Gamma$. 
Thus $d(e, g \cdot e) = d(e, g s_{14}) = 5$ implies $g s_{14}$ has word length 5 with respect to the generators 
$s_{12}$, $s_{23}$, $s_{34}$, $s_{13}$, $s_{24}$ in $J_4'$, and $\pi(g s_{14})$ must be an odd permutation. 
Since $\pi(s_{14}) = (14)(23)$, it implies that $\pi(g)$ is an odd permutation, contradicting $\pi(g) = e$. 

As a result, $D(e)$ contains the interior of $\widetilde{D}_1$. 
It concludes that $\widetilde{D}_1 = \mathrm{int}(\widetilde{D} \cap \mathcal{G})$ is the Dirichlet domain with center $e$ on $\mathcal{G}$ with respect to the action $\Gamma_1$ of $PJ_4$.
\end{proof}

In the following, the side of $\partial \widetilde{D}$ with endpoints $v, v'$ is denoted by $s(v,v')$. 

\begin{proof}[Proof of Theorem~\ref{PJ4 presentation}]
By Lemma~\ref{D tilde is Dirichlet domain}, it follows that the interior of $\widetilde{D}$ is the Dirichlet domain on $\mathcal{C}$, and thus $\widetilde{D}$ is a Dirichlet polygon on $\mathcal{C}$ with respect to the action $\Gamma$ of $PJ_4$ on $\mathcal{C}$. 
This is because the action $\Gamma$ on $\mathcal{C}$ is a natural extension of $\Gamma_1$ on $\mathcal{G} = \mathcal{C}^{(1)}$, which is an isometric action on $\mathcal{C}$, and $\widetilde{D}$ is the convex hull in $\mathcal{C}$ of $\widetilde{D} \cap \mathcal{G}$.

Since a group acting freely and cocompactly on a topological space is discontinuous, together with Proposition~\ref{prop42}, it follows that $PJ_4$ is discontinuous. 
By Proposition~\ref{C4_32_isomH^2}, $\mathcal{C}$ is geodesically connected, geodesically complete, and finitely compact. 
Thus, applying Proposition~\ref{Dirichlet domain to fundamental domain} to the interior of $\widetilde{D}$, we obtain that $\widetilde{D}$ is a locally finite fundamental polygon for $PJ_4$. By Proposition~\ref{D generates G}, $PJ_4$ is generated by the elements $g_{1}, \dots, g_{10}$ and their inverses.

Next, we will obtain the relations for $g_1, \dots, g_{10}$ by applying Proposition~\ref{poincare's thm}.  
The generators $g_1, \dots, g_{10}$ of $PJ_4$ define identifications on the sides of the polygon $\widetilde{D}$.  
For example, by
\begin{align*}
    g_{1}(s_{34}s_{12})&=a_{1}s_{14}(s_{34}s_{12})s_{14} \\
    &=(s_{13}s_{24}s_{12}s_{34})s_{14}s_{14}(s_{12}s_{34}) \\
    &=s_{13}s_{24}s_{12}s_{34}(s_{12}s_{34}) \\
    &=s_{13}s_{24}s_{12}s_{34}(s_{34}s_{12}) \\
    &=s_{13}s_{24},
\end{align*}
\begin{align*}
    g_{1}(s_{34}s_{13})&=a_{1}s_{14}(s_{34}s_{13})s_{14} \\
    &=(s_{13}s_{24}s_{12}s_{34})s_{14}s_{14}(s_{12}s_{24}) \\
    &=s_{13}s_{24}s_{12}s_{34}(s_{12}s_{24}) \\
    &=s_{13}s_{24}s_{12}(s_{12}s_{34}s_{24}) \\
    &=s_{13}s_{24}(s_{34}s_{24}) \\
    &=s_{13}s_{24}(s_{24}s_{23}) \\
    &=s_{13}s_{23},
\end{align*}
we see that 
\[
g_1 : s(s_{34} s_{12} , s_{34} s_{13} ) \mapsto s( s_{13} s_{24} , s_{13} s_{23} ).
\]

%The above equalities imply that $g_{1}$, as an isometry on $\mathcal{C}$, identifies the sides %$\overline{v(s_{34}s_{12})v(s_{34}s_{13}})$ 
%$s(s_{34}s_{12},s_{34}s_{13})$ and %$\overline{v(s_{13}s_{24})v(s_{13}s_{23})}$. 
%$s( s_{13}s_{24} ,s_{13}s_{23} )$. 

In the same way, we have the following. 
\[
\begin{array}{rll}
g_{2} &: 
s ( s_{24}s_{13} , s_{24}s_{13}s_{23} ) 
&\mapsto 
s( s_{13}s_{24} , s_{13}s_{24}s_{23} )\\
g_{3} &:
s( s_{34}s_{23} , s_{34}s_{23}s_{13}) 
&\mapsto 
s ( s_{13}s_{34} , s_{13}s_{34}s_{24} )\\
g_{4} &:
s( s_{24}s_{34}) , s_{24}s_{12} ) 
&\mapsto 
s( s_{13}s_{34} , s_{13}s_{12} )\\
g_{5} &:
s( s_{13}s_{23} , s_{12}s_{23} )
&\mapsto 
s( s_{23}s_{12} , s_{23}s_{13} ) \\
g_{6} &:
s ( s_{34}s_{13} , s_{34}s_{13}s_{12} ) 
& \mapsto 
s( s_{23}s_{12} , s_{23}s_{12}s_{34} ) \\
g_{7} &:
s ( s_{12}s_{24} , s_{12}s_{24}s_{34} ) 
& \mapsto 
s( s_{23}s_{34} , s_{23}s_{34}s_{12} ) \\
g_{8} &: 
s( s_{24}s_{23} , s_{34}s_{23} ) 
& \mapsto 
s( s_{23}s_{34} , s_{24}s_{34} ) \\
g_{9} &: 
s( s_{12}s_{23} , s_{12}s_{23}s_{24} ) 
& \mapsto 
s( s_{24}s_{12} , s_{24}s_{12}s_{13} ) \\
g_{10} &: 
s( s_{12}s_{34} , s_{12}s_{24} ) 
& \mapsto 
s( s_{24}s_{13} , s_{24}s_{23} )
\end{array}
% \begin{array}{ll}
% g_{2}(s_{24}s_{13}) = s_{13}s_{24} \qquad & 
% g_{2}(s_{24}s_{13}s_{23}) = s_{13}s_{24}s_{23}\\
% g_{3}(s_{34}s_{23}) = s_{13}s_{34} & 
% g_{3}(s_{34}s_{23}s_{13}) = s_{13}s_{34}s_{24} \\
% g_{4}(s_{24}s_{34}) = s_{13}s_{34} & 
% g_{4}(s_{24}s_{12}) = s_{13}s_{12} \\
% g_{5}(s_{13}s_{23}) = s_{23}s_{12} & 
% g_{5}(s_{12}s_{23}) = s_{23}s_{13} \\
% g_{6}(s_{34}s_{13}) = s_{23}s_{12} & 
% g_{6}(s_{34}s_{13}s_{12}) = s_{23}s_{12}s_{34} \\
% g_{7}(s_{12}s_{24}) = s_{23}s_{34} & 
% g_{7}(s_{12}s_{24}s_{34}) = s_{23}s_{34}s_{12} \\
% g_{8}(s_{24}s_{23}) = s_{23}s_{34} & 
% g_{8}(s_{34}s_{23}) = s_{24}s_{34} \\
% g_{9}(s_{12}s_{23}) = s_{24}s_{12} & 
% g_{9}(s_{12}s_{23}s_{24}) = s_{24}s_{12}s_{13} \\
% g_{10}(s_{12}s_{34}) = s_{24}s_{13} & 
% g_{10}(s_{12}s_{24}) = s_{24}s_{23}
% \end{array}
\]

It follows that $g_1$, \dots $g_{10}$ give identifications on $\widetilde{D}$. 
The identified polygon $\widetilde{D}^*$ under the identifications $g_1$, \dots, $g_{10}$ is a closed surface. 
%Furthermore, the following discussion concludes $\widetilde{D}$ is complete and satisfies the cycle condition, and so $\widetilde{D}$ is Poincar\'{e} polygon. 

In the following, we will confirm that $\widetilde{D}$ is complete and satisfies the cycle condition, and so $\widetilde{D}$ is a Poincar\'{e} polygon. 

Let $u_1$ be the vertex $s_{13}s_{24}s_{23}$ and $s_1$ the side 
$s( s_{13}s_{24}s_{23} , s_{13}s_{24} )$. 
Then $g_{2}^{-1}$ maps $u_1$ to the vertices 
$u_2:=s_{24}s_{13}s_{23}$ 
and $s_1$ to the side 
$s( s_{24}s_{13}s_{23} , s_{24}s_{13} )$. 
Let 
%a vertex $u_{2}$ as $g_{2}^{-1}(u_1)=s_{24}s_{13}s_{23}$ and 
$s_2$ be the side 
$s ( s_{24}s_{13}s_{23} , s_{24}s_{12})$ 
%$\overline{v(s_{24}s_{13}s_{23})v(s_{24}s_{12})}$, 
which is the edge sharing $u_2$ with $g_2^{-1}(s_1)$. 
Then, by $g_{9}^{-1}$, the vertex 
$u_{2}=s_{24}s_{13}s_{23} =s_{24}s_{12}s_{13}$ 
and the side $s_2$ are mapped to the vertex 
$u_3 : = s_{12}s_{23}s_{24} = s_{12}s_{24}s_{34}$ 
and $s( s_{12}s_{23}s_{24} , s_{12}s_{23} ) = s ( s_{12}s_{24}s_{34} , s_{12}s_{23} )$. 
%$\overline{v(s_{12}s_{24}s_{34})v(s_{12}s_{23})}$. 
%We denote the vertex $g_{9}^{-1}(u_2)$ by $u_3$ and 
Let $s_3 := s( s_{12}s_{24}s_{34} , s_{12} s_{24} )$ be the side sharing $u_3$ with $g_{9}^{-1}(s_2)$. 
The generator $g_7$ maps the vertex $u_3$ and the side $s_3$ to 
$u_4 := s_{23}s_{34}s_{12}$ and
$s ( s_{23}s_{34}s_{12} , s_{23}s_{34} )$. 
%\overline{v(s_{23}s_{34}s_{12})v(s_{23}s_{34})}$. 
%We denote the vertex $g_7(s_3)$ by $u_4$ and 
Let $s_4$ be the side 
$s ( s_{23}s_{34}s_{12} , s_{23}s_{12}) =s ( s_{23}s_{12}s_{34} , s_{23}s_{12})$ 
which is the edge sharing $u_4$ with $g_7(s_3)$. 
By the generator $g_6^{-1}$, the vertex $u_4$ and the side $s_4$ are mapped to 
$u_5 := s_{34}s_{13}s_{12}$ and 
$s ( s_{34}s_{13}s_{12} , s_{34}s_{13})$. 
%We denote the vertex $g_{6}^{-1}(u_4)$ by $u_5$ and 
Let $s_5$ be the side 
$s( s_{34}s_{13}s_{12} , s_{34}s_{23}) 
= s(s_{34}s_{23}s_{13} ,  s_{34}s_{23} )$, 
which is the side sharing $u_5$ with $g_6^{-1}(s_4)$. 
The generator $g_3$ maps the vertex $u_5$ to $u_1$, and $s_5$ to the side 
$ s ( s_{13}s_{34}s_{24} , s_{13}s_{34}  ) =
s ( s_{13}s_{24}s_{23} , s_{13}s_{34})$, 
which is the edge on $\partial \widetilde{D}$ sharing $u_1$ with $s_1$.

Therefore, the sequence of generators 
\[
\{ g_2^{-1} , g_9^{-1} , g_7 , g_6^{-1} , g_3 \}
\]
is a cycle of generators and the sequence of vertices $\{u_1, \dots , u_5\}$ is a cycle of vertices. 
Moreover, all of the angles $\alpha(u_i)$, made by the two sides around $u_i$ and measured from inside $\widetilde{D}$, are equal to $\frac{2\pi}{5}$. 
Thus, we obtain $\sum_{i=e}^{5} \alpha (u_i) = 2 \pi$. 

In the same way, we obtain the other cycles as follows. 

\[
\{ g_{4}^{-1} , g_{8}^{-1} , g_3 \} , 
\qquad 
\{ s_{13}s_{34} , s_{24}s_{34} , s_{34}s_{23} \}
\]
The angles $\alpha( s_{13}s_{34} )$ and $\alpha(s_{34}s_{23})$ are equal to $\frac{3\pi}{5}$ and the angle $\alpha( s_{24}s_{34} )$ is $\frac{4\pi}{5}$, and so, the sum of the angles is equal to $2 \pi$.

\[
\{ g_{4}^{-1} , g_9^{-1} , g_5 \} , 
\qquad 
\{ s_{13}s_{12} , s_{24}s_{12} , s_{12}s_{23} \}
\]
The angle $\alpha( s_{13}s_{12} )$ is $\frac{4\pi}{5}$ and the angles $\alpha( s_{24}s_{12} )$ and $\alpha(s_{12}s_{23})$ are equal to $\frac{3\pi}{5}$, and so, the sum of the angles is equal to $2 \pi$. 

\[
\{ g_{6}^{-1} , g_1 , g_5 \} , 
\qquad 
\{ s_{23}s_{12} , s_{34}s_{13} , s_{13}s_{23} \}
\]
The angles $\alpha( s_{23}s_{12} )$ and $\alpha(s_{34}s_{13})$ are equal to $\frac{3\pi}{5}$ and the angle $\alpha( s_{13}s_{23} )$ is $\frac{4\pi}{5}$, and so, the sum of the angles is equal to $2 \pi$. 

\[
\{ g_{7}^{-1} , g_{10} , g_8 \} , 
\qquad 
\{ s_{23}s_{34} , s_{12}s_{24} , s_{24}s_{23} \}
\]
The angles $\alpha( s_{23}s_{34} )$ and $\alpha(s_{12}s_{24})$ are equal to $\frac{3\pi}{5}$ and the angle $\alpha( s_{24}s_{23} )$ is $\frac{4\pi}{5}$, and so, the sum of the angles is equal to $2 \pi$. 

\[
\{ g_2 , g_1^{-1} , g_{10} \} , 
\qquad 
\{ s_{24}s_{13} , s_{13}s_{24} , s_{34}s_{12} \}
\]
The angles $\alpha( s_{24}s_{13} )$ and $\alpha(s_{13}s_{24})$ are equal to $\frac{3\pi}{5}$ and the angle $\alpha( s_{34}s_{12} )$ is $\frac{4\pi}{5}$, and so, the sum of the angles is equal to $2 \pi$. 

It concludes that $\widetilde{D}$ satisfies the cycle condition, and so $\widetilde{D}$ is a Poincar\'{e} polygon. 
Applying Proposition~\ref{poincare's thm}, the compositions of generators $g_{3} g_6^{-1} g_{7} g_{9}^{-1} g_{2}^{-1}$, $g_{3} g_{8}^{-1} g_{4}^{-1}$, $g_{5} g_{9}^{-1} g_{4}^{-1}$, $g_{5} g_{1} g_{6}^{-1}$, $g_8 g_{10} g_7^{-1}$, $g_{10} g_{1}^{-1} g_{2}$ form a complete set of relations. 
     As a result, we obtain the presentation 
\[
      \left\langle 
      g_1, \cdots, g_{10} \left| 
      \begin{array}{l}
      g_{3} g_6^{-1} g_{7} g_{9}^{-1} g_{2}^{-1} 
      = g_{3} g_{8}^{-1} g_{4}^{-1}
      = g_{5} g_{9}^{-1} g_{4}^{-1} \\
      = g_{5} g_{1} g_{6}^{-1}
      = g_8 g_{10} g_7^{-1}
      = g_{10} g_{1}^{-1} g_{2} =e
% g_1g_{10}^{-1}g_2^{-1}=g_9g_5^{-1}g_4=g_5g_1g_6^{-1}\\
% =g_8g_{10}g_7^{-1}=g_8g_3^{-1}g_4=g_{2}g_{9}g_{7}^{-1}g_{6}g_{3}^{-1} =e
      \end{array}
     \right.
     \right\rangle
\]

Moreover, by substituting 
$g_1=g_2 g_{10}$, 
$g_5=g_4 g_9$, 
$g_6=g_5 g_1$, 
$g_7=g_8 g_{10}$, 
$g_3=g_4 g_8$ to 
$g_{3} g_6^{-1} g_{7} g_{9}^{-1} g_{2}^{-1} =e$, 
we obtain the following presentation. 
\begin{align*}
\langle g_2, g_4, g_8, g_9, g_{10} \mid g_{2}g_{9}g_{10}^{-1}g_{8}^{-1}g_{4}g_{9}g_{2}g_{10}g_{8}^{-1}g_{4}^{-1} =e
\rangle
\end{align*}
\end{proof}

We remark that, by calculating the Euler characteristic, we see that the surface obtained from $\widetilde{D}$ by gluing its sides via $g_j$'s is homeomorphic to the connected sum of the five projective plane.

\section{Proof of corollaries}

In this section, we give proofs of corollaries. 

\begin{proof}[Proof of Corollary~\ref{BCL}]
Set 
\[
\begin{array}{l}
G := \langle \alpha, \beta, \gamma, \delta, \epsilon \mid \alpha \gamma \epsilon \beta \epsilon \alpha^{-1} \delta^{-1} \beta \gamma \delta^{-1} =e \rangle , \\
G':= \langle g_2, g_4, g_8, g_9, g_{10} \mid g_{2}g_{9}g_{10}^{-1}g_{8}^{-1}g_{4}g_{9}g_{2}g_{10}g_{8}^{-1}g_{4}^{-1} =e \rangle .
\end{array}
\]

Consider the following correspondence between the generators of $G$ and $G'$. 
\[
\begin{array}{ll}
   \alpha \longmapsto g_4 g_9 \\
   \beta \longmapsto g_2 g_{10} \\
   \gamma \longmapsto g_{10}^{-1} \\
   \delta \longmapsto g_4 \\
   \epsilon \longmapsto g_8^{-1} g_4 g_9 
\end{array}
\]

Then the following calculation shows that the correspondence implies a well-defined homomorphism, say $f: G \to G'$. 
  \begin{align*}
      \alpha \gamma \epsilon \beta \epsilon \alpha^{-1} \delta^{-1} \beta \gamma \delta^{-1} 
      \longmapsto \ & g_4 g_9 g_{10}^{-1} g_{8}^{-1} g_4 g_9 g_2 g_{10} g_8^{-1} g_4 g_9 (g_4 g_9)^{-1} g_4^{-1} g_2 g_{10} g_{10}^{-1} g_4^{-1}  \\
      &= g_9 g_{10}^{-1} g_{8}^{-1} g_4 g_9 g_2 g_{10} g_8^{-1} g_4^{-1} g_2 
      = e 
  \end{align*}

This map $f$ is a surjection since the following hold. 
\begin{align*}
f(\beta \gamma) &= g_2 g_{10} g_{10}^{-1} = g_2 \\
f(\delta) &= g_4 \\
f(\alpha \epsilon^{-1}) &= g_4 g_9 g_9^{-1} g_4^{-1} g_8 = g_8\\
f(\delta^{-1} \alpha) &= g_4^{-1} g_4 g_9 = g_9 \\
f(\gamma^{-1}) &= g_{10}
\end{align*}

In the same way, we can construct the inverse of $f$. 
The correspondence 
\[
\begin{array}{ll}
   g_2 \longmapsto \beta \gamma \\
   g_4 \longmapsto \delta \\
   g_8 \longmapsto \alpha \epsilon^{-1} \\
   g_9 \longmapsto \delta^{-1} \alpha \\
   g_{10} \longmapsto \gamma^{-1} 
\end{array}
\]
induces a well-defined homomorphism $g : G' \to G$ by the following. 
 \begin{align*}
g_{10}^{-1} g_{8}^{-1} g_4 g_9 g_2 g_{10} g_8^{-1} g_4^{-1} g_2 g_9 
\longmapsto \ &
\gamma \epsilon \alpha^{-1} \delta \delta^{-1} \alpha \beta \gamma \gamma^{-1} \epsilon \alpha^{-1} \delta^{-1} \beta \gamma \delta^{-1} \alpha \\
&= \gamma \epsilon \beta \epsilon \alpha^{-1} \delta^{-1} \beta \gamma \delta^{-1} \alpha 
= e
\end{align*}

%Since the generators $\alpha$, \dots, $\epsilon$ satisfies the following, 
Then it follows that $g \circ f = id$ from the following. 
 \[
 \begin{array}{ll}
g \circ f (\alpha) &= g(g_4 g_9) = \delta \delta^{-1} \alpha = \alpha \\
g \circ f (\beta) &= g(g_2 g_{10}) = \beta \gamma \gamma^{-1} = \beta \\
g \circ f (\gamma) &= g(g_{10}^{-1}) = (\gamma^{-1})^{-1} = \gamma \\
g \circ f (\delta) &= g(g_4) = \delta \\
g \circ f (\epsilon) &= g(g_8^{-1} g_4 g_9) = \epsilon \alpha^{-1} \delta \delta^{-1} \alpha = \epsilon 
 \end{array}
 \]
 As a result, $f$ is shown to be an isomorphism. 
\end{proof}

\begin{proof}[Proof of Corollary~\ref{Pi_1_of_ConnSum}]
Set 
\[
\begin{array}{ll}
G & := \langle \alpha_1 , \alpha_2 , \alpha_3 , \alpha_4 , \alpha_5 \mid 
\alpha_1^2\alpha_2^2\alpha_3^2\alpha_4^2\alpha_5^2 = e \rangle , \\[12pt]
G'& := 
\left\langle 
g_1, \cdots, g_{10} \left| 
\begin{array}{l}
g_{3} g_6^{-1} g_{7} g_{9}^{-1} g_{2}^{-1} 
= g_{3} g_{8}^{-1} g_{4}^{-1}
= g_{5} g_{9}^{-1} g_{4}^{-1} \\
= g_{5} g_{1} g_{6}^{-1}
= g_8 g_{10} g_7^{-1}
= g_{10} g_{1}^{-1} g_{2} =e
\end{array}
\right.
\right\rangle 
% \\
% & =
% \langle g_1, \cdots, g_{10} \mid 
% \begin{array}{l}
% g_1 g_{10}^{-1}g_2^{-1}
% = g_9g_5^{-1}g_4
% = g_5g_1g_6^{-1} \\
% = g_8g_{10}g_7^{-1}
% = g_8g_3^{-1}g_4
% = g_{2}g_{9}g_{7}^{-1}g_{6}g_{3}^{-1} 
% \end{array}
% \rangle .
\end{array}
\]

Consider the following correspondence. 
\[
\begin{array}{ll}
    \alpha_1 \longmapsto g_1^{-1} = g_{10}^{-1} g_2^{-1} \\
    \alpha_2 \longmapsto g_2 g_{10} g_5^{-1} g_8 g_3^{-1} = g_2 g_{10} g_9^{-1} g_4^{-2} \\
    \alpha_3 \longmapsto g_4 \\
    \alpha_4 \longmapsto g_9 g_{10}^{-1} \\
    \alpha_5 \longmapsto g_8^{-1} g_6 = g_8^{-1} g_4 g_9 g_2 g_{10} 
\end{array}
\]
Then, by the next calculation, we have a well-defined homomorphism $f : G \to G'$. 
    \begin{align*}
\alpha_1^2 \alpha_2^2 \alpha_3^2 \alpha_4^2 \alpha_5^2 
\longmapsto \ & 
g_1^{-2} ( g_2 g_{10} g_5^{-1} g_8 g_3^{-1} )^2 g_4^2 (g_9 g_{10}^{-1})^2 (g_8^{-1} g_6)^2 \\ 
&= g_1^{-2} (g_1 g_5^{-1} g_4^{-1})^2 g_4^2 (g_9 g_{10}^{-1})^2 (g_8^{-1} g_6)^2 \\ 
&= g_1^{-1} g_5^{-1} g_4^{-1} g_1 (g_5^{-1} g_4) g_9 g_{10}^{-1} g_9 (g_{10}^{-1} g_8^{-1}) g_6 g_8^{-1} g_6 \\ 
&= g_1^{-1} g_5^{-1} g_4^{-1} g_1 ((g_9^{-1}) g_9) g_{10}^{-1} g_9 (g_7^{-1}) g_6 g_8^{-1} g_6 \\ 
&=  (g_1^{-1} g_5^{-1}) g_4^{-1} (g_1 g_{10}^{-1}) g_9 g_7^{-1} g_6 g_8^{-1} g_6 \\ 
&= g_6^{-1} g_4^{-1} g_3 (g_6^{-1} g_6) g_8^{-1} g_6 \\ 
&= g_6^{-1} g_4^{-1} g_4 g_8 g_8^{-1} g_6 = e 
\end{align*}

This map $f$ is a surjection from the following. 
\begin{align*}
f(\alpha_2 \alpha_3^2 \alpha_4) 
&= g_2 g_{10} g_9^{-1} g_4^{-2} g_4^2 g_9 g_{10}^{-1} 
=  g_2 \\
f(\alpha_3) &= g_4 \\
f(\alpha_3^{-1} \alpha_2^{-1} \alpha_1^{-2} \alpha_5^{-1} ) 
&= g_4^{-1} g_4^2 g_9 g_{10}^{-1} g_2^{-1} (g_2 g_{10})^2 g_{10}^{-1} g_2^{-1} g_9^{-1} g_4^{-1} g_8 
= g_8\\
f(\alpha_3^{-2} \alpha_2^{-1} \alpha_1^{-1}) 
&= g_4^{-2} g_4^2 g_9 g_{10}^{-1} g_2^{-1} g_2 g_{10} 
= g_9 \\
f((\alpha_2 \alpha_3^2 \alpha_4)^{-1} \alpha_1^{-1}) 
&= g_2^{-1} g_2 g_{10} 
= g_{10}
\end{align*}

On the other hand, the correspondence 
\[
\begin{array}{lll}
g_2 & \longmapsto  \alpha_2 \alpha_3^2 \alpha_4 \\
g_4  & \longmapsto \alpha_3 \\
g_8 & \longmapsto \alpha_3^{-1} \alpha_2^{-1} \alpha_1^{-2} \alpha_5^{-1}\\
g_9 & \longmapsto \alpha_3^{-2} \alpha_2^{-1} \alpha_1^{-1} \\
g_{10} &\longmapsto \alpha_4^{-1} \alpha_3^{-2} \alpha_2^{-1} \alpha_1^{-1}
\end{array}
\]
induces a well-defined homomorphism $g : G' \to G$ by the following. 
\[
\begin{array}{ll}
& g_{2} g_{9} g_{10}^{-1} g_{8}^{-1} g_{4} g_{9} g_{2} g_{10} g_{8}^{-1} g_{4}^{-1} \\
& \qquad \qquad 
\longmapsto 
(\alpha_2 \alpha_3^2 \alpha_4) (\alpha_3^{-2} \alpha_2^{-1} \alpha_1^{-1}) (\alpha_1 \alpha_2 \alpha_3^2 \alpha_4) (\alpha_5 \alpha_1^2 \alpha_2 \alpha_3) \alpha_3 (\alpha_3^{-2} \alpha_2^{-1} \alpha_1^{-1}) \\
& \qquad \qquad\qquad   
(\alpha_2 \alpha_3^{2} \alpha_4 ) 
(\alpha_4^{-1} \alpha_3^{-2}) \alpha_2^{-1} \alpha_1^{-1} 
(\alpha_5 \alpha_1^2 \alpha_2 \alpha_3 ) \alpha_3^{-1} \\
& \qquad \qquad \qquad = \alpha_2 \alpha_3^2 \alpha_4 \alpha_4 \alpha_5 \alpha_1^2 (\alpha_2 \alpha_3 \alpha_3 \alpha_3^{-2} \alpha_2^{-1}) \alpha_1^{-1} (\alpha_2 \alpha_2^{-1}) \alpha_1^{-1} \alpha_5 \alpha_1^2 \alpha_2 \\
& \qquad \qquad \qquad = \alpha_2 \alpha_3^2 \alpha_4 \alpha_4 \alpha_5 (\alpha_1^2 \alpha_1^{-1} \alpha_1^{-1}) \alpha_5 \alpha_1^2 \alpha_2 \\
& \qquad \qquad \qquad 
= \alpha_2 \alpha_3^2 \alpha_4 \alpha_4 \alpha_5 \alpha_5 \alpha_1^2 \alpha_2 = e
%&  \qquad = \alpha_2 \alpha_3^2 \alpha_4^2 \alpha_5^2 \alpha_1^2 \alpha_2 
% = e 
\end{array}
\]
Then it follows that $g \circ f = id$ as; 
  \[
  \begin{array}{ll}
        g \circ f (\alpha_1) = g (g_{10}^{-1} g_{2}^{-1}) 
        = \alpha_1 ( \alpha_2 \alpha_3^2 \alpha_4 \alpha_4^{-1} \alpha_3^{-2} \alpha_2^{-1} ) 
        = \alpha_1 \\
        g \circ f (\alpha_2) = g (g_2 g_{10} g_9^{-1} g_4^{-2}) 
        = \alpha_2 (\alpha_3^2 \alpha_4 \alpha_4^{-1} \alpha_3^{-2}) (\alpha_2^{-1} \alpha_1^{-1} \alpha_1 \alpha_2) (\alpha_3^2 \alpha_3^{-2}) 
        = \alpha_2 \\
        g \circ f (\alpha_3) = g (g_4) 
        = \alpha_3 \\
        g \circ f (\alpha_4) = g (g_9 g_{10}^{-1}) 
        = \alpha_3^{-2} \alpha_2^{-1} \alpha_1^{-1} \alpha_1 \alpha_2 \alpha_3^2 = \alpha_4 \\
        g \circ f (\alpha_5) = g ( g_8^{-1} g_4 g_9 g_2 g_{10} ) 
        = \alpha_5 (\alpha_1^2 \alpha_2 \alpha_3 \alpha_3 \alpha_3^{-2} \alpha_2^{-1} \alpha_1^{-1}) \alpha_2 (\alpha_3^2 \alpha_4 \alpha_4^{-1} \alpha_3^{-2}) \alpha_2^{-1} \alpha_1^{-1} \\
        \qquad \qquad \,
        = \alpha_5 (\alpha_1 \alpha_2 \alpha_2^{-1} \alpha_1^{-1}) 
        = \alpha_5 
  \end{array}
  \]  
This concludes that $f$ gives an isomorphism. 
\end{proof}

\begin{figure}[htb]
    \centering
\begin{overpic}[width=.75\textwidth]{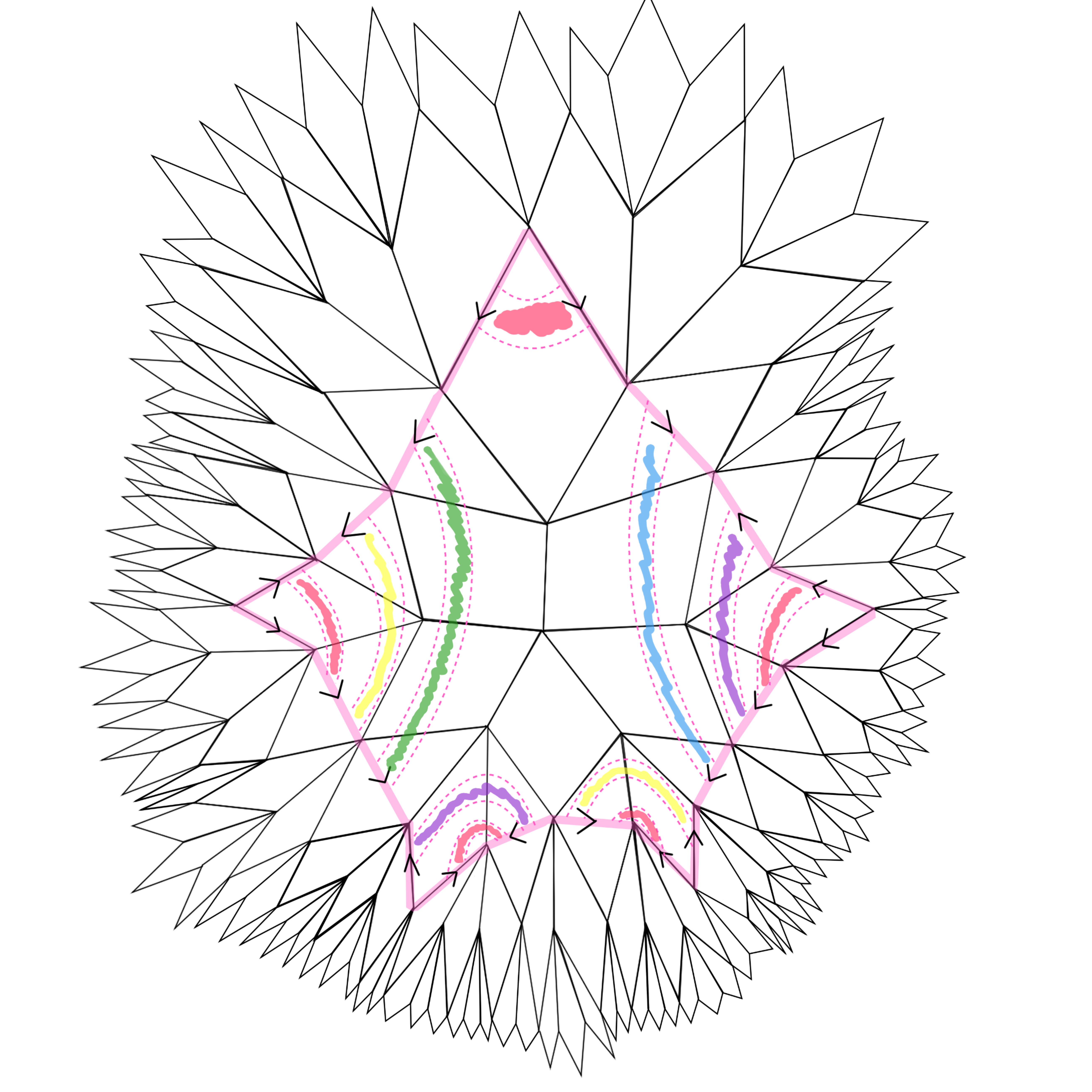}
\end{overpic}
\caption{}
\label{Colored Mobius band}
\end{figure}

\begin{figure}[htb]
    \centering
\begin{overpic}[width=.75\textwidth]{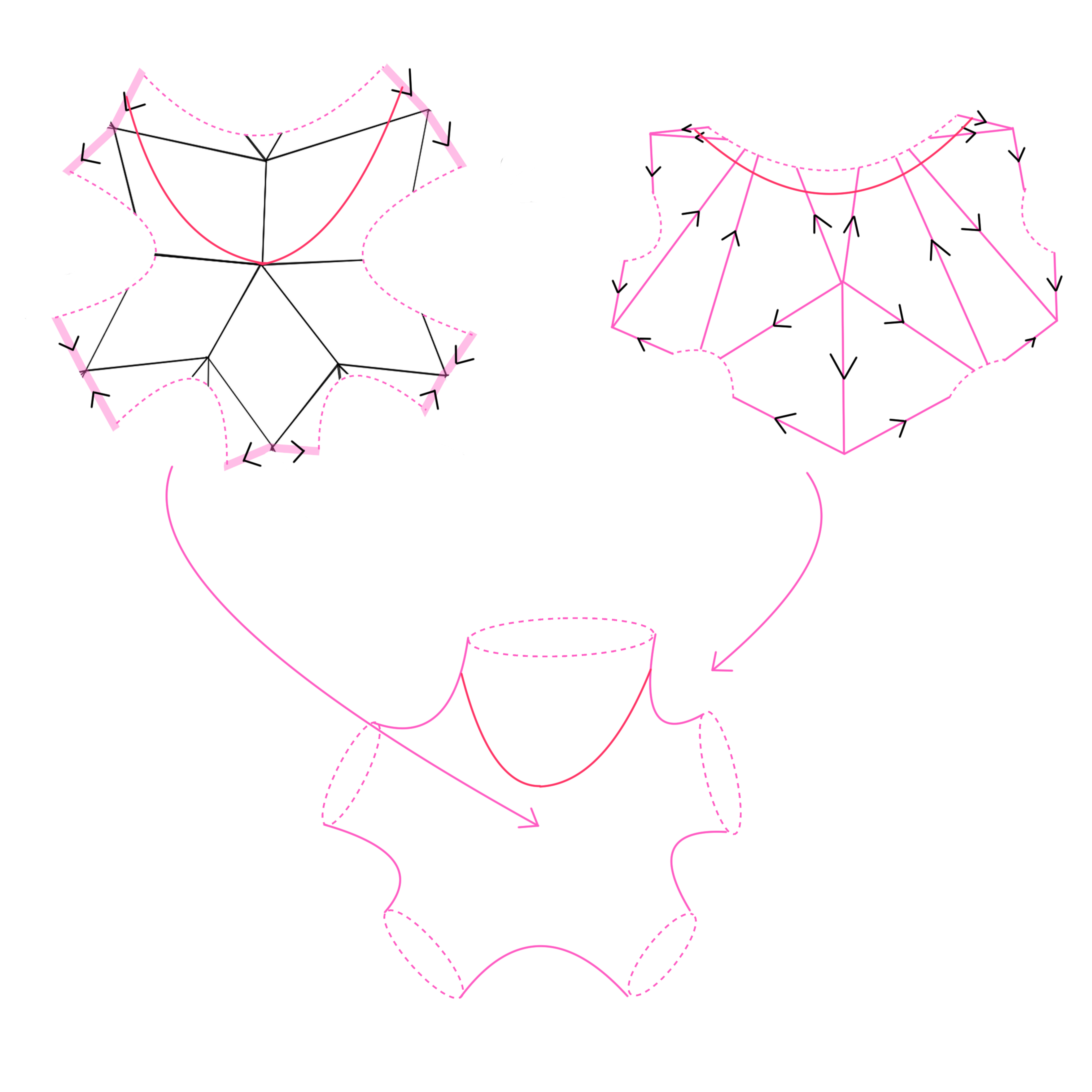}
\put(50,30){\textcolor{red}{$\alpha_2$}}
\end{overpic}
\caption{}
\label{C_4^23/PJ_4}
\end{figure}

\begin{remark} 
By virtue of geometric descriptions, we could find the correspondence between $\alpha_i$ and $g_j$ in the proof above as follows.  
The painted regions in Figure~\ref{Colored Mobius band} represent five M\"{o}bius bands embedded in $\mathcal{C}/PJ_4$.  
After removing them, the remaining parts can be reassembled to form a sphere with five holes.  
See Figure~\ref{C_4^23/PJ_4}.  
Then, we can confirm that $\mathcal{C} / PJ_4$ is homeomorphic to the connected sum $\#_5 \mathbb{RP}^2$ of five real projective planes.  
By setting the closed curve $\alpha_2$ in $\#_5 \mathbb{RP}^2$ shown in Figure~\ref{C_4^23/PJ_4}, we obtain a correspondence  
\[
    \alpha_2 \longmapsto g_8^{-1} g_4 g_9 g_2 g_{10} \in PJ_4
\]
by considering a lift of $\alpha_2$ in $\mathcal{C}$.  
Similarly, we obtain the corresponding elements in $PJ_4$ for $\alpha_1$, $\alpha_3$, $\alpha_4$, and $\alpha_5$.  
\end{remark}

% \begin{figure}[htb]
%     \centering
% \begin{overpic}[width=.6\textwidth]{How_to_find_2.pdf}
% \end{overpic}
% \caption{}
% \label{soko_futa}
% \end{figure}

\section*{Acknowledgments}
We would like to thank Takuya Sakasai for his valuable feedback and Carl-Fredrik Nyberg Brodda for inspiring future research themes. 

\bibliographystyle{amsplain}
\bibliography{ref}

\end{document}